\RequirePackage[ngerman,english]{babel}
\documentclass[12pt,a4paper]{amsart}

\usepackage{amsfonts, amsmath, amssymb, amsthm}

\usepackage[ngerman,english]{babel}
\usepackage{ngerman}
\usepackage{setspace}

\theoremstyle{plain}
\newtheorem{thm}{Theorem}[section]
\newtheorem*{thm*}{Theorem}
\newtheorem{lem}[thm]{Lemma}
\newtheorem{prop}[thm]{Proposition}
\newtheorem{cor}[thm]{Corollary}

\theoremstyle{definition}

\newtheorem{exmp}[thm]{Example}

\theoremstyle{remark}
\newtheorem{rem}[thm]{Remark}
\newtheorem{qst}[thm]{Question}

\usepackage{algorithmic,algorithm}

\usepackage{txfonts}
\usepackage{amssymb}
\usepackage{mathbbol}
\usepackage{ae}
\usepackage{bbm}
\usepackage[mathscr]{eucal} % euler fonts for calligraphic letters -> \mathscr{}
\usepackage{xspace}
\usepackage{url}

\usepackage{color}
\usepackage{graphicx}
\usepackage{overpic}
\usepackage{subfigure}
\usepackage{textcomp}
\usepackage{booktabs}
\usepackage{paralist}
\usepackage[a4paper,scale=0.732, marginratio={1:1, 9:10}, ignoreall]{geometry}

\numberwithin{equation}{section}
\numberwithin{figure}{section}

\DeclareMathOperator{\Vertices}{Vert}

\newcommand\qm[1]{``#1''}

\newcommand\cut{\cap}
\newcommand\Cut{\bigcap}
\newcommand\union{\cup}
\newcommand\Union{\bigcup}

\newcommand\isom{\cong}

\newcommand\cA{{\mathcal A}}
\newcommand\cB{{\mathcal B}}
\newcommand\cF{{\mathcal F}}
\newcommand\cH{{\mathcal H}}

\newcommand\cM{{\mathcal M}}

\newcommand\NN{{\mathbb N}}
\newcommand\RR{{\mathbb R}}

\newcommand\SetOf[2]{\left\{#1\vphantom{#2}\,\right.\left|\,\vphantom{#1}#2\right\}}
\newcommand\smallSetOf[2]{\{#1\,|\,#2\}}

\newcommand\aff{\operatorname{aff}}
\newcommand\conv{\operatorname{conv}}

\newcommand\vol{\operatorname{vol}}

\newcommand\vertt{\Vertices}
\newcommand\intt{\operatorname{int}}
\newcommand\relint{\operatorname{relint}}

\newcommand\card[1]{\left|#1\right|}

\newcommand{\Gr}{{\mathrm{Gr}}}

\newcommand{\Dr}{{\mathrm{Dr}}}
\providecommand{\nonnegRR}{{\RR_{\ge0}}}

\DeclareMathOperator{\pos}{pos}
\providecommand\cprime{$'$}
\providecommand{\subdiv}{{\Sigma}}
\providecommand{\subdivG}{{\tilde\Sigma}}

\providecommand{\envelope}[2]{{{\mathscr{E}}_{#1}(#2)}}
\providecommand{\lift}[2]{{{\mathscr{L}}_{#1}(#2)}}
\providecommand{\tightspan}[2]{{{\mathscr{T}}_{#1}(#2)}}

\providecommand{\subdivision}[2]{{\subdiv_{#1}(#2)}}
\providecommand{\subdivisionG}[2]{{\subdivG_{#1}(#2)}}
\providecommand{\Asubdivision}[2]{{\subdiv_{#1}(#2)}}

\providecommand{\Hypersimplex}[2]{{\Delta(#1,#2)}}

\providecommand{\transpose}[1]{{#1}^{\sf T}}
\providecommand{\scp}[2]{\langle{#1},{#2}\rangle}

\providecommand{\SecondaryPolytope}[1]{{\rm SecPoly}(#1)} % use $\Sigma$ for individual subdivision
\providecommand{\SplitPoly}[1]{{\rm SplitPoly}(#1)}

\renewcommand{\phi}{\varphi}

\providecommand{\ChamberComplex}[1]{{\rm Chamber}(#1)}

\usepackage[ansinew]{inputenc}
\usepackage[T1]{fontenc}

\graphicspath{ {fig/} }

\title{On the Facets of the Secondary Polytope}
\address{Fachbereich Mathematik, TU Darmstadt, 64289 Darmstadt, Germany}
\email{sherrmann@mathematik.tu-darmstadt.de}

\author{Sven Herrmann}

\begin{document}

\begin{abstract}
The secondary polytope of a point configuration~$\cA$ is a polytope whose face poset is isomorphic to the poset of all regular subdivisions of~$\cA$. While the vertices of the secondary polytope -- corresponding to the triangulations of~$\cA$ -- are very well studied, there is not much known about the facets of the secondary polytope.

The splits of a polytope, subdivisions with exactly two maximal faces, are the simplest examples of such facets and the first that were systematically investigated. The present paper can be seen as a continuation of these studies and as a starting point of an examination of the subdivisions corresponding to the facets of the secondary polytope in general. As a special case, the notion of $k$"=split is introduced as a possibility to classify polytopes in accordance to the complexity of the facets of their secondary polytopes. An application to matroid subdivisions of hypersimplices and tropical geometry is given.
\end{abstract}
\maketitle

\section{Introduction}
A subdivision of a point configuration~$\cA$ is a collection $\subdiv$ of subsets of~$\cA$ (the \emph{faces} of $\subdiv$) such that the union of the convex hull of all of the faces equals the convex hull of~$\cA$ and such that the intersection of two faces of $\subdiv$ is a face of both. Subdivisions and especially triangulations (i.e., subdivisions into simplices) occur in various parts of mathematics; for an overview see the first chapter of the monograph~\cite{Triangulations} by De Loera, Rambau, and Santos. One way to construct polytopal subdivisions of~$\cA$ is the following: Let $w:\cA \to \RR$ be a function assigning a weight to each element of~$\cA$. By lifting each $a\in\cA$ according to its weight and projecting the lower faces of the resulting polytope down to $\conv \cA$, one obtains a subdivision of~$\cA$. Such subdivisions are called \emph{regular}. It is an important structural result by Gel{\cprime}fand, Kapranov, and Zelevinsky \cite{MR1073208} (see also \cite[Chapter~7]{MR1264417}) that there exists a polytope $\SecondaryPolytope \cA$, called the \emph{secondary polytope} of $\cA$, whose vertices are in bijection with the regular triangulations of~$\cA$. Moreover, they showed that the face poset of $\SecondaryPolytope \cA$ is isomorphic to the poset of all regular subdivisions of~$\cA$ ordered by refinement. In this way, the facets of $\SecondaryPolytope \cA$ correspond to those regular subdivisions of~$\cA$ that can only be coarsened by the trivial subdivision. The aim of this paper is to start an investigation of those coarsest subdivisions.

In \cite{MR2502496} Joswig and the author studied the notion of \emph{split} of a polytope~$P$, generalizing earlier work on finite metric spaces by Bandelt and Dress~\cite{BandeltDress86}; see also Hirai~\cite{MR2252108}. These are the simplest possible (non"=trivial) subdivisions of a polytope one can think of, namely those with exactly two maximal faces. These splits are special kinds of the facets of the secondary polytope of~$P$. We will generalize these ideas in two ways: First, we will study splits of general point configurations that do not have to be in general position; almost all the results about splits generalize trivially to this more general case. The second generalization is more interesting: We will study a much bigger class of facets of the secondary polytope of a point configuration~$\cA$: the \emph{$k$"=splits}. The $k$"=splits of~$\cA$ are those coarsest subdivisions of~$\cA$ that have exactly one interior face of codimension $k-1$. One of our main results is the assertion that all of these subdivisions are indeed regular, hence facets of the secondary polytope.

As a next step,  we will study general coarsest subdivisions of point configurations that are not necessarily $k$"=splits. In doing so, we will use the notion of \emph{tight span} of a polyhedral subdivision, which was also introduced in \cite{MR2502496} and which originates in the theory of finite metric spaces \cite{MR753872,Isb}. The tight span of a subdivision $\subdiv$ is the polyhedral complex dual to the interior faces of $\subdiv$. This concept allows us to investigate how complicated coarsest subdivisions with a given number~$k$ of maximal faces can get, and we give a classification of all corresponding tight spans for small~$k$.

One case where one is much more interested in the facets of the secondary polytopes rather than the vertices, is the study of matroid subdivisions. It was shown by Speyer~\cite{MR2448909} that the space of all regular matroid subdivisions of the hypersimplex $\Hypersimplex kn$ is the space of all $(k-1)$"=dimensional tropical linear spaces in tropical $(n-1)$"=dimensional space. This space is (a close relative of) the tropical analogue of the Grassmannian; see~\cite{SpeyerSturmfels04}. Since triangulations can never be matroid subdivisions, one key step in the study of all matroid subdivisions is the determination of the coarsest matroid subdivisions, which generate the space of all such subdivisions. We will show that $3$"=splits of $\Hypersimplex kn$ are matroid subdivisions for all $l\leq k$.

This paper is organized as follows. In the beginning, we give basic definitions and results about subdivisions and point configurations used in the sequel including the generalization of the theory of tight spans from polytopes to point configurations. In Section~\ref{sec:pc-polys}, we will give two results justifying that subdivisions of point configurations are -- in principle -- not more complicated than subdivisions of polytopes. First, all secondary polytopes arising for point configurations arise for polytopes, too, second, any tight span occurring for a point configuration occurs for some polytope. In the end of the section, we will show that each polytope can be the tight span of some subdivision of another polytope. In Section~\ref{sec:k-splits}, we will give our examples and results about general $k$"=splits and the fifth section investigates the tight spans of $k$"=subdivisions, general coarsest subdivisions with $k$ maximal faces. After some general discussions of the possible tight spans for $k$"=subdivisions, we give classifications of the tight spans of $k$"=subdivisions for small $k$ and show that not all polytopes can be the tight span of some $k$"=subdivision. After the discussion of $3$"=splits of hypersimplices and their matroid subdivision, we conclude the paper with a list of open questions.

The author would like to thank the anonymous referees for their helpful comments and suggestions.

\section{Subdivisions of Point Configurations}\label{sec:subdivs}

A \emph{point configuration} is a finite multiset~$\cA\subset\RR^d$. By a \emph{multiset} we mean a collection whose members may appear multiple times. Throughout, we suppose that~$\cA$ has the maximal dimension $d$, where the dimension of a point configuration~$\cA$ is defined as the dimension of the affine hull $\aff \cA$. A \emph{subdivision} $\subdiv$ of~$\cA$ is a collection of subconfigurations of~$\cA$ satisfying the following three conditions (see~\cite[Section~2.3]{Triangulations}):
\begin{itemize}
\item(SD1)\label{sd:containement} If $F\in\subdiv$ and $\bar F$ is a face of $F$, then $\bar F\in \subdiv$.
\item(SD2)\label{sd:union} $\conv \cA =\Union_{F\in \subdiv} \conv F$.
\item(SD3)\label{sd:intersection} If $F, \bar F\in\subdiv$, then $\relint(\conv F) \cut \relint( \conv \bar F) =\emptyset$.
\end{itemize}

A point configuration $F\subset \cA$ is called a \emph{face} of~$\cA$ if there exists a supporting hyperplane $H$ of $\conv \cA$ such that $F=\cA\cut H$.

Given a subdivision $\subdiv$ of a point configuration~$\cA$, we can look at the polyhedral complex $\subdivG:=\smallSetOf{\conv F}{F\in\subdiv}$. This is a polyhedral subdivision of the polytope $\conv \cA$ possibly with additional vertices. Note that for two different subdivisions $\subdiv\not=\subdiv'$ of a point configuration~$\cA$, we can have $\subdivG=\subdivG'$; see Example~\ref{ex:hexagon:pc}. We will sometimes call~$\subdivG$ a \emph{geometric subdivision} of~$\cA$ in order to distinguish it from the subdivision~$\subdiv$.

If $P$ is a polytope, we can consider the point configuration~$\cA(P):=\vertt{P}$ consisting of the vertices of $P$. A \emph{subdivision} $\subdiv$ of $P$ is defined as a subdivision of~$\cA(P)$. This implies that all points used in $\subdiv$ are vertices of $P$. Furthermore, for subdivisions $\subdiv$, $\subdiv'$ of $P$,  $\subdivG=\subdivG'$ is equivalent to $\subdiv=\subdiv'$, so we do not have to distinguish between $\subdiv$ and the geometric subdivision $\subdivG$ for polytopes. %In fact, our main interest is in the study of subdivisions of polytopes, however

\subsection{Regular Subdivisions and Tight Spans}

Given a weight function $w:\cA \to \RR$ we consider the lifted polyhedron
\begin{equation*}
    \lift{w}{\cA} \ := \conv\SetOf{(w(a),a)}{a\in\cA} \, + \, \nonnegRR (1,0,\dots,0) \subset \RR\times \RR^d\, .
  \end{equation*}
  
The \emph{regular subdivision} $\Asubdivision w \cA$ of~$\cA$ with respect to $w$ is obtained by taking the sets $\smallSetOf{b\in\cA}{(w(b),b)\in F}$ for all lower faces $F$ (with respect to the first coordinate; by definition, these are exactly the bounded faces) of $\lift w \cA$. So the elements of $\subdivisionG w \cA$ are the projections of the bounded faces of $\lift w \cA$ to the last $d$ coordinates.

Furthermore, we define the \emph{envelope} of~$\cA$ with respect to $w$ as

\[
\envelope{w}{\cA} \ := \ \SetOf{x\in\RR\times\RR^{d}}{\scp {(1,a)} x \ge -w \text{ for all } a\in\cA}
\]

and the \emph{tight span} $\tightspan w \cA$ of~$\cA$ as the complex of bounded faces of $\envelope{w}{\cA}$. From this, one derives that for two lifting functions $w_1,w_2$ we have that $\tightspan {w_1} \cA=\tightspan {w_2} \cA$ implies $\subdivisionG {w_1} \cA=\subdivisionG {w_2} \cA$ but not necessarily $\Asubdivision {w_1} \cA=\Asubdivision {w_2} \cA$; see Example~\ref{ex:square}, also for illustrations of the concepts of envelope and tight span.

The following proposition, which is a direct generalization of \cite[Proposition~2.3]{MR2502496} and can be shown in the same way, gives the relation between tight spans and regular subdivisions.

\begin{prop}\label{prop:duality}
  The polyhedron $\envelope{w}{\cA}$ is affinely equivalent to the polar dual of the polyhedron $\lift w \cA$.
  Moreover, the face poset of $\tightspan{w}{\cA}$ is anti"=isomorphic to the face poset of
  the interior lower faces (with respect to the first coordinate) of $\lift{w}{\cA}$.
\end{prop}

So the (inclusion) maximal faces of the tight span $\tightspan{w}{\cA}$ correspond to the (inclusion) minimal interior faces of $\Asubdivision w \cA$. Here, a face of $\Asubdivision w \cA$ is an \emph{interior face} if it is not entirely contained in the boundary of $\conv \cA$. We will be especially interested in those subdivisions that have exactly one minimal interior face; we say that these subdivisions have the \emph{G"=property}. By Proposition~\ref{prop:duality}, a subdivision has the G"=property if and only if its tight span is (the complex of faces of) a single polytope. Furthermore, we will say that a point configuration~$\cA$ has the \emph{G"=property} if all coarsest subdivisions of~$\cA$, that is, those subdivisions that cannot be refined non"=trivially, have the G"=property.

\begin{rem}
  The G"=property is related to the notion of Gorenstein polytopes \cite{MR2275581}, Gorenstein simplicial complexes, and Gorenstein rings \cite{MR1251956,MR1453579} as follows. A simplicial complex~$\Delta$ is \emph{Gorenstein} if the polynomial ring $K[\Delta]$ is a Gorenstein ring. It was shown by Joswig and Kulas \cite[Proposition~24]{JoswigKulas08} that a regular triangulation~$\subdiv$ (considered as simplicial complex) is Gorenstein if and only if its tight span has a unique maximal cell, that is, if and only if $\subdiv$ has the G"=property. By a result of Bruns and Römer \cite[Corollary~8]{MR2275581}, a polytope (satisfying some additional properties) is Gorenstein if and only if it has some Gorenstein triangulation. So if we have such a polytope $P$ and a triangulation of $P$ with the G"=property, then $P$ is Gorenstein. It would be interesting to explore how general subdivisions with the G"=property and polytopes with the G"=property translate into the commutative algebra setting of Gorenstein simplicial complexes and Gorenstein rings.

\end{rem}

We call a sum $w_1+w_2$ of two weight functions of a point configuration~$\cA$ \emph{coherent} if
\begin{equation}\label{eq:coherent}
  \envelope{w_1}{\cA}+\envelope{w_2}{\cA} \ = \ \envelope{w_1+w_2}{\cA} \, .
\end{equation}
(Note that $\subseteq$ holds for any weight functions.)
We get the following corollary of Proposition~\ref{prop:duality} translating this property into the language of regular subdivisions.

\begin{cor}\label{cor:coherent2}
 A decomposition $w=w_1+w_2$ of weight functions of~$\cA$ is coherent if and only if
  the subdivisions
  $\Asubdivision{w_1}{\cA}$ and $\Asubdivision{w_2}{\cA}$ have a common refinement.
\end{cor}

We postpone the proof of Corollary~\ref{cor:coherent2} to the end of this section since it uses the theory of secondary polytopes, which we will discuss in the next subsection. Before this, though, we would like to mention that also most of the other elementary results proved in \cite[Section~2]{MR2502496} are true for general point configurations, too. However, sometimes one has to be careful whether one has to consider $\Asubdivision w \cA$ or $\subdivisionG w\cA$.

\begin{figure}[htb]
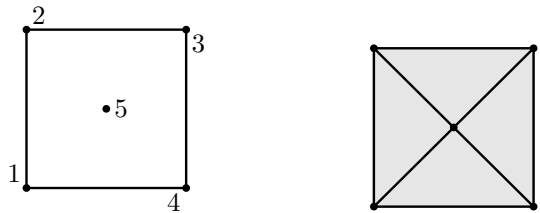
\centering
  \includegraphics{square.0}
  \hskip 2cm
  \includegraphics{square.1}
  \caption{The point configuration of Example~\ref{ex:square}.}
  \label{fig:square-pc}
\end{figure}

\begin{exmp}\label{ex:square}
We consider the point configuration~$\cA$ whose elements are the columns of the matrix
\[
V \ = \
  \begin{pmatrix}
%    1 & 1 & 1 & 1 & 1\\
    0 & 0 & 2 & 2 & 1\\
    0 & 2 & 0 & 2 & 1
  \end{pmatrix}
  \]
consisting of the vertices of a square together with its center (see Figure~\ref{fig:square-pc}) and the weight functions $w_1=(1,0,0,0,0)$, $w_2=(0,0,1,0,0),$ $\bar{w_1}=(1,0,0,0,1)$, and $\bar{w_2}=(0,0,1,0,1)$. A computation shows that the envelope of  $w_1$ and $\bar w_1$ is a three"=dimensional unbounded polyhedron with two vertices and four rays: 
\[
\envelope{w_1}{\cA} =\conv\{(0,0,0),(-1,1/2,1/2)\}+\pos\{(2,-1,0),(2,0,-1),(0,0,1),(0,1,0)\}\,.
\]
(Here $\pos S$ is the set of all positive linear combinations of a given set $S$.) So $\tightspan{w_1}{\cA} = \tightspan{\bar{w_1}}{\cA}$ is the polyhedral complex consisting of the line segment $[(0,0,0),(-1,1/2,1/2)]$, its two vertices, and the empty set. Similarly, $\tightspan{w_2}{\cA} = \tightspan{\bar{w_2}}{\cA}$ is the face poset of $[(0,0,0),(0,-1/2,-1/2)]$.
We now have that $\subdivisionG {w_i} \cA=\subdivisionG {\bar {w_i}} \cA$, but $\Asubdivision {w_i} \cA\not=\Asubdivision {\bar {w_i}} \cA$ for $i\in\{1,2\}$. The geometric subdivisions $\subdivisionG {\bar {w_1}} \cA$ and $\subdivisionG {\bar {w_2}} \cA$ have a common refinement, the subdivision depicted in Figure~\ref{fig:square-pc} on the right, just as  $\subdivisionG {{w_1}} \cA$ and $\subdivisionG {{w_2}} \cA$. The corresponding subdivision is also the common refinement of $\Asubdivision {w_1} \cA$ and $\Asubdivision {w_2} \cA$, but $\Asubdivision {\bar{w_1}} \cA$ and $\Asubdivision {\bar{w_2}} \cA$ do not have a common refinement. This agrees with the fact that  $w_1+w_2$ is coherent, whereas $\bar{w_1}+\bar{w_2}$ is not, and verifies Corollary~\ref{cor:coherent2} in this case.
\end{exmp}

 \begin{exmp}\label{ex:hexagon:pc}
Consider the point configuration $\cH$ forming a hexagon with an interior point, consisting of the columns of the matrix \[
  V \ = \
  \begin{pmatrix}
 %   1 & 1 & 1 & 1 & 1 & 1 & 1\\
    0 & 1 & 2 & 2 & 1 & 0 & 1\\
    0 & 0 & 1 & 2 & 2 & 1 & 1
  \end{pmatrix}\,,
  \]
   and the weight functions $w=(0,0,1,1,0,0,0)$ and $\bar w=(0,0,1,1,0,0,1)$. A direct computation shows that $$\tightspan{w}\cH\ =\ \tightspan{\bar w}\cH\ =\ \tightspan{w_1}{H}\ =\ [0,(1,-1,0)]\;.$$
The geometric subdivisions $\subdivisionG {w} \cH$ and $\subdivisionG {\bar w} \cH$ agree, but the subdivisions $\Asubdivision {w} \cH$ and $\Asubdivision {\bar w} \cH$ are not equal: The former has the maximal faces $\{1,2,5,6,7\}$ and $\{2,3,4,5,7\}$, but the latter the maximal faces $\{1,2,5,6\}$ and $\{2,3,4,5\}$. Here, the numbers correspond to the columns of the matrix $V$. So $\Asubdivision {\bar w} \cH$ is strictly finer than $\Asubdivision {w} \cH$. %One can see this also be generalizing the  \emph{coherency index} \cite[Equation~(2)]{MR2502496} to point configurations: One computes $\alpha^w_{\bar w}=0$ and $\alpha_w^{\bar w}=1$. This shows that  and gives us a counterexample to Corollary~\ref{cor:equal-subdivision} for the case of geometric subdivisions.

%The decomposition $(w,\bar w)$ is coherent and $\Asubdivision {w+\bar w} \cH=\Asubdivision {\bar w} \cH$, $\subdivisionG {w+\bar w} \cH=\subdivisionG {\bar w} \cH$ in accordance with Corollary~\ref{cor:coherent}.
\end{exmp}

\begin{rem}\label{rem:non-regular}
%\begin{enumerate}
%\item As we defined $k$"=subdivisions so far, they are abstract subdivisions. However, we will also call the geometric subdivision $\Sigma(\Delta)$ a $k$"=subdivision if $\Delta$ is a $k$"=subdivision. But note that there might be another abstract subdivision $\Delta'$ that is not a $k$"=subdivision with $\Sigma(\Delta')=\Sigma(\Delta)$. Since we are mainly interested in the tight spans of the subdivisions in this section, we will normally work with a geometric subdivisions $\Sigma$, but the results are true as well for $\Delta_\cA(\Sigma)$.
%\item 
So far, we only defined the tight span for regular subdivisions. However, for any subdivision $\subdiv$ of a point configuration~$\cA$ one can define the tight span $\tightspan{\subdiv}{\cA}$ as the abstract polyhedral complex that is dual to the complex of interior faces of~$\subdiv$. For regular subdivisions, the usual tight span is a realization of this abstract polyhedral complex by Proposition~\ref{prop:duality}.% When we proof general properties of tight spans makes perfect sense also for non"=regular subdivisions.
%\end{enumerate}
\end{rem}

\subsection{Secondary Polytopes}\label{sec:secondary-polytope}

%The set of all regular subdivisions of a point configuration~$\cA$ is known to have an
%interesting structure (see \cite[Chapter 5]{Triangulations} for the details): For $w: \cA\to \RR$ consider the set $S[w]\subset \RR^n$ of all weight functions that define the same regular subdivision of~$\cA$. This set is called the \emph{secondary cone} of~$\cA$ with respect to $w$.  It can be shown (for instance, see \cite[Corollary 5.2.10]{Triangulations}) that $S[w]$ is indeed a polyhedral cone and that the set of all $S[w]$ (for all $w$) forms a polyhedral fan $\SecondaryFan{\cA}$, called the \emph{secondary fan} of~$\cA$. Since an affine linear function defines the trivial subdivision on each point configuration~$\cA$, all cones $S[w]$ have a lineality space of dimension $d+1$ and $\SecondaryFan \cA$ can be regarded in the quotient space  isomorphic to $\RR^{n-d-1}$.

The \emph{secondary polytope} of a point configuration~$\cA$ was first defined by Gel{\cprime}fand, Kapranov, and Zelevinsky. They showed \cite[Theorem 1.7]{MR1073208} that there exists a polytope, the \emph{secondary polytope}
$\SecondaryPolytope{\cA}$ of~$\cA$, whose face poset is isomorphic to the poset of all regular subdivisions of~$\cA$. This polytope admits a realization as the convex hull of
the so"=called \emph{GKZ"=vectors} of all triangulations of~$\cA$. The GKZ"=vector
$x_\subdiv\in \RR^{\cA}$ of a triangulation $\subdiv$ of $\cA$ is defined as ${(x_\subdiv)}_a:=\sum_{S}
\vol S$ for all $a\in\cA$, where the sum ranges over all full"=dimensional
simplices $S\in\subdiv$ that contain~$a$. A dual description of the secondary polytope in terms of its facets was given by Lee \cite[Section 17.6, Result~4]{MR1730170}. Each facet"=defining inequality is obtained explicitly from a weight function of the corresponding coarsest regular subdivision.

%Most of the results about tight spans and 2"=splits of polytopes remain valid for point configurations. However, in this more general context, we do not only have to deal with 2"=splits but also with $1$"=splits of point configurations. A \emph{$1$"=split} of a point configuration~$\cA$ is a subdivision of~$\cA$ with exactly one maximal face containing all but one of the points in~$\cA$. It will turn out that these one"=splits are needed to prove a generalization of the 2"=split Decomposition Theorem for point configurations.

The normal fan of $\SecondaryPolytope \cA$ is called the \emph{secondary fan} of~$\cA$. Actually, an (open) cone in the secondary fan is given by the set of all weight functions that define the same regular subdivision of~$\cA$; see, for example, \cite[Chapter~5]{Triangulations} for a detailed discussion of secondary fans (and secondary polytopes). 

There is a nice way to construct the secondary fan of a point configuration given by Billera, Filliman, and Sturmfels~\cite[Section~4]{MR1074022}. We describe this construction very briefly and refer to~\cite{MR1074022,BilleraSturmfels93} for the details. The key ingredient for this construction is the \emph{Gale transform} of a point configuration; see \cite[Section~5.4]{MR1976856} or \cite[Chapter~6]{MR1311028}. Let~$\cA$ be a point configuration, $n:=\card \cA$, and $V$ the $n\times(d+1)$"=matrix whose rows are the points $(1,a)$ for all $a\in\cA$. Now consider an $n\times(n-d-1)$"=matrix $V^\star$ of full rank $n-d-1$ satisfying
$\transpose{V}V^\star=0$; that is, the columns of $V^\star$ form a basis of the kernel of
$\transpose{V}$.  Then the rows of $V^\star$ form a vector configuration~$\cB$ in $\RR^{n-d-1}$. This configuration is called the \emph{Gale dual} or \emph{Gale transform} of~$\cA$. The multiset~$\cB$ has the same number of elements as~$\cA$ and if $a\in\cA$ corresponds to the $i$th row of $V$, the element of~$\cB$ corresponding to the $i$th row of~$V^\star$ is called $a^\star$. Note that~$\cB$ may be a proper multiset even if~$\cA$ does not have any multiple points.

The \emph{chamber complex} $\ChamberComplex \cB$ of~$\cB$ is the coarsest polyhedral complex that covers $\pos \cB=\RR^{n-d-1}$ and that refines all triangulations of~$\cB$. Details and a combinatorial study of the chamber complex can be found in \cite{AGZ}; see also \cite[Section~5.3]{Triangulations}. The relation of the chamber complex of the Gale dual~$\cB$ with the secondary polytope of~$\cA$ is the following.

\begin{thm}[{\cite[Theorem~3.1]{BilleraSturmfels93}}]\label{thm:secondaryDuality}
  The chamber complex $\ChamberComplex{\cB}$ is anti"=isomorphic to the boundary complex of
  the secondary polytope $\SecondaryPolytope{\cA}$.
\end{thm}

%Any Gale dual of~$\cA$ is uniquely determined up to affine equivalence.  Each vector $v\in V$ corresponds to a row vector $v^\star$ of $V^\star$, called the \emph{vector dual} to $v$.  Throughout we will assume that all dual vectors are either zero or have unit Euclidean length.  If $v^\star$ is zero then all vectors other than $v$ span a linear hyperplane not containing $v$. We call $V$ \emph{proper} if $V^\star$ does not contain any zero vectors. In the primal view, this means that $\conv V$ is not a pyramid.  For the remainder of this section we will assume that $V$ is proper whence $V^\star$ can be identified with a configuration of $n$ points on the unit sphere $\Sph^{n-d-2}$.  Notice that these $n$ points are not necessarily pairwise distinct. Repetitions may occur even if the vectors in $V$ are pairwise distinct.

This bijection can be made explicit as follows: Let $w:\cA \to \RR$ be a weight function. This weight~$w$ is identified with the vector $\tilde w:=\sum_{a\in\cA} w(a)\cdot a^\star\in \RR^{n"=d"=1}$. The regular subdivision $\Asubdivision w \cA$ is uniquely determined by $\tilde w$ since one can show that for two weight functions $w_1,w_2$ with $\tilde w_1=\tilde w_2$ one has that $w_1-w_2$ is an affine linear function, which obviously induces the trivial subdivision on~$\cA$. The regular subdivision $\Asubdivision w \cA$ can now be determined from $\tilde w$ and  $\ChamberComplex \cB$: A subconfiguration $F\subset \cA$ is an element of $\Asubdivision w \cA$ if and only if $\tilde w \in \intt \pos \smallSetOf{a^\star}{a\not\in F}$; see \cite[Lemma~3.2]{BilleraSturmfels93}.

\begin{proof}[Proof of Corollary~\ref{cor:coherent2}]
By \cite[Corollary~2.4]{MR2502496}, a decomposition $w=w_1+w_2$ of weight functions for a polytope $P$ is coherent if and only if the subdivision $\subdivision{w}{P}$ is the common refinement of the subdivisions $\subdivision{w_1}{P}$ and $\subdivision{w_2}{P}$. The proof of this statement can be literally generalized to point configurations. So it remains to prove that the existence of a common refinement of $\Asubdivision{w_1}{\cA}$ and $\Asubdivision{w_2}{\cA}$ implies the coherence. In terms of the secondary polytope, the existence of a refinement of $\Asubdivision{w_1}{\cA}$ and $\Asubdivision{w_2}{\cA}$ implies that the intersection of the corresponding faces of $\SecondaryPolytope \cA$ is non"=empty. So, by Theorem~\ref{thm:secondaryDuality}, the chambers of $\ChamberComplex \cB$ with $\tilde w_1, \tilde w_2$ in their relative interior lie in a common chamber $C$. However, the chamber~$C$ then also contains $\tilde w_1 + \tilde w_2$, which can be retranslated to the statement that $\subdivision{w}{\cA}$ is the  common refinement of $\subdivision{w_1}{\cA}$ and $\subdivision{w_2}{\cA}$. Hence $w=w_1+w_2$ is coherent.
\end{proof}

\section{Point Configurations and Polytopes}\label{sec:pc-polys}

In this section, we will give two results concerning the \qm{complexity} of subdivisions of point configurations relative to polytopes. Both results say that, in principle, subdivisions of point configurations do not get more complicated than those of polytopes.

The first result is a corollary of Theorem~\ref{thm:secondaryDuality}; compare, for example, \cite[Theorem~4.2.35]{Triangulations}. We include a simple proof of this statement.

\begin{thm}\label{thm:pc-poly-secondary}
Let~$\cA$ be a $d$"=dimensional point configuration with $n$ points. Then there exists a $d+m$"=dimensional polytope $P$ with $n+m$ vertices such that~$\cA$ and~$\cA(P)$ have isomorphic secondary polytopes and $m\leq n$.
\end{thm}

\begin{proof}
By \cite[Section~5.4, Theorem~2]{MR1976856}, a vector configuration~$\cB$ is the Gale dual of a polytope if and only if every open halfspace whose boundary contains the origin contains at least two elements of~$\cB$. Let~$\cB$ be the Gale dual of~$\cA$. Since~$\cB$ is positively spanning, every such open halfspace contains at least one element of~$\cB$. If there exists some halfspace with exactly one element, say $b$, we add a copy of $b$ to~$\cB$, and we repeat this step until we have two elements in each halfspace. The derived vector configuration~$\cB'$ is the Gale dual of some polytope $P$. However, we have $\ChamberComplex \cB=\ChamberComplex{\cB'}$ by the definition of the chamber complex, hence Theorem~\ref{thm:secondaryDuality} shows that the secondary polytopes of $\cA$ and $P$ are isomorphic. Since the number of points added is at most $n$, we also get the proposed bound.
\end{proof}

\begin{rem}
\begin{enumerate}
\item That the bound proposed in Theorem~\ref{thm:pc-poly-secondary} is sharp, can be seen by the following trivial example. Let $\cA$ be the $0$"=dimensional point configuration consisting of $n$ copies of a single point. Then the Gale dual $\cB$ of $\cA$ consists of $n$ linear independent vectors and we have to add a copy for each of them. The resulting polytope ~$P$ is the $n$"=dimensional cross polytope.
\item The polytope~$P$ constructed from the point configuration~$\cA$ in the proof of Theorem~\ref{thm:pc-poly-secondary} is a multiple one"=point suspension of~$\cA$; see \cite[Section~4.2.5]{Triangulations}.
\item In the same manner as in the proof of Theorem~\ref{thm:pc-poly-secondary}, starting with the Gale dual $\cB$ of any point configuration~$\cA$ one arrives at point configurations~$\cA'$ with isomorphic secondary polytopes. This shows that for any point configuration there exist infinitely many (non"=isomorphic) proper point configurations that have the same secondary polytope. (A point configuration is \emph{proper} if it is not a pyramid; a point configuration $\cA$ and the pyramid over $\cA$ obviously have isomorphic secondary polytopes.)
\item If $\ChamberComplex \cB$ contains a ray $r$ with $r\not=\pos b$ for all $b\in \cB$, then one can add any $c$ with $r=\pos c$ to $\cB$ without changing the chamber complex. The existence of such a ray in $\ChamberComplex \cB$ is equivalent to the existence of a regular coarsest subdivision $\Sigma$ of~$\cA$ that does not contain~$\cA \setminus \{a\}$ as a maximal cell for some $a\in \cA$. In particular, this condition is satisfied if~$\cA$ has more coarsest subdivisions than elements. Hence, in this case, there exist (finitely many) point configurations with the same secondary polytope as~$\cA$ that are not obtained via one"=point suspensions.
\item The construction in the proof of Theorem~\ref{thm:pc-poly-secondary} is related to the Lawrence construction; see Billera and Munson~\cite[Section~2]{MR782051}. However, rather than adding the negative of an existing vector in the Gale dual, we add a copy.
\end{enumerate}
\end{rem}

The second result is of a very different nature. Whereas Theorem~\ref{thm:pc-poly-secondary} talks about the structure of the collection of all regular subdivisions of~$\cA$ and $P$, it does not give any information about the relation between the individual subdivisions of~$\cA$ and $P$. For example, the number of maximal cells of the subdivision usually changes. The following result, however,  concerns the combinatorics of an individual subdivision of~$\cA$ in terms of its tight span: By considering point configurations instead of polytopes one does not allow more possibilities for the tight spans.

\begin{prop}\label{prop:pc-to-poly}
Let~$\cA\subset \RR^{d}$ be a point configuration with $n$ points and $\Asubdivision{w}{\cA}$ a regular subdivision of~$\cA$. Then there exists a polytope $P\subset \RR^{d+1}$ with $2n$ vertices together with a regular subdivision~$\Asubdivision{w'}{P}$ of~$P$ such that $\tightspan{w'}{P}$ is affinely isomorphic to $\tightspan w \cA$. Furthermore, if $\Asubdivision{w}{\cA}$ is a coarsest subdivision of~$\cA$, then $\Asubdivision{w'}{P}$ is a coarsest subdivision of~$P$.
\end{prop}

\begin{proof}
By possibly deleting some points from~$\cA$, we can assume that~$\cA$ does not have any multiple points and that for each cell $F\in\subdivision{w}{\cA}$ all $a\in F$ are vertices of $\conv F$. Furthermore, we assume that $w<0$. Then we define the polytope $P\subset\RR^{d+1}=\RR^{d}\times\RR$ as
\[
P\ :=\ \conv\SetOf{(a,\pm w(a))}{a\in\cA}\,.
\]
From our assumption that every $a\in\cA$ is the vertex of some $F\in\subdivision{w}{\cA}$, it follows that all lifted points $(w(a),a)$ are vertices of $\lift{w}{\cA}$; and so from $w<0$ it follows that all points $(a,\pm w(a))$ are vertices of $P$. We define a weight function $w':\cA(P) \to \RR$ as $w'(a,\pm w(a))=w(a)$. From the definition of the envelope, we directly get that $x\in\envelope{w}{\cA}$ implies $(x,0)\in\envelope{w'}{P}$ and that $(x,x')\in\envelope{w'}{P}$ implies $x\in\envelope{w}{\cA}$. We will now show that $\tightspan{w'}P=\tightspan w \cA \times \{0\}$, which implies the claim.

Since the vertices of $\envelope w \cA$ are the vertices of $\tightspan w \cA$, it suffices to show that $(v,v')\in\RR^{d+1}\times \RR$ is a vertex of $\envelope{w'}P$ if and only if $v'=0$ and $v$ is a vertex of $\envelope w \cA$. So let first $v$ be a vertex of $\envelope w \cA$. Then there exists a $(d+1)$"=element set~$C\subset \cA$ such that $v$ is the unique solution $x\in \RR^{d+1}$ of the linear system $\scp {(1,a)} x = -w(a)$ for all $a\in\cB$. This implies that $(v,0)$ is the unique solution $(x,x')\in \RR^{d+1}\times \RR$ to the system $\scp {(1,a)} x \pm w(a)x'=-w(a)$ for all $a\in C$, and so $(v,0)$ is a vertex of $\envelope{w'}P$.

On the other hand, consider a vertex $(v,v')\in\RR^{d+1}\times (\RR\setminus \{0\})$ of
\[
\envelope{w'}P=\SetOf{(x,x')\in\RR^{d+1}\times \RR}{\scp {(1,a)} x \pm w(a)\geq w(a)x' \text{ for all } a\in\cA}\,.
\]
 Suppose that there exists some $p,q\in\cA$ with
\begin{align}
\scp {(1,p)} v + w(p)v'\ &=\ -w(p)\quad \text{and}\label{eq:pm:1}\\
\scp {(1,q)} v - w(q)v'\ &=\ -w(q)\,.\label{eq:pm:2}
\end{align}
 Since $v\in\envelope{w}{\cA}$ we have $\scp {(1,p)} v \geq -w(p)$ and $\scp {(1,q)} v \geq -w(q)$. Furthermore, by our assumption, we have $w(p),w(q)<0$. So Equation~\eqref{eq:pm:1} yields $v'> 0$, and Equation~\eqref{eq:pm:2} yields $v'<0$, a contradiction. So we can assume that we only have equality in \qm{$+$}"=inequalities. Hence, we find a $(d+2)$"=element set~$C\subset\cA$ such that $(v,v')$ is the unique solution $(x,x')$ of the linear system $\scp{(1,a)} x + w(a)x'=-w(a)$ for all $a\in\cB$. However, a solution to this system is $(0,-1)$, which is not an element of~$\envelope{w'}{P}$ (since it does not fulfill any of the \qm{$-$}"=inequalities). This contradiction finishes the proof of the first assertion. 
 
What remains to show is that the subdivision $\subdivision{w'}{P}$ cannot be refined non"=trivially if this was the case for $\subdivision w \cA$. Suppose there exists some non"=trivial coarsening~$\subdiv'$ of~$\subdivision{w'}{P}$. It is easily checked that $\subdiv':=\smallSetOf{C\cut(\RR^{d+1}\times\{0\})}{C\in\subdiv'}$ is a subdivision of~$\cA\times\{0\}$, so we would also have a subdivision of~$\cA$ that coarsens $\subdivision w \cA$ non"=trivially. %This contradicts the assumption that $\subdivision w \cA$ was a $k$"=subdivision.
\end{proof}

%\begin{rem}
%One can show that
%\[
%\envelope{w'}{P}\ =\ \SetOf{(x,x')\in \RR^{d+1}\times \RR}{x\in\envelope w \cA\text{ and }|x'|\leq-\max_{a\in\cA}\frac{\scp a x}{w(a)}-1}.
%\]
%\end{rem}

In contrast to Theorem~\ref{thm:pc-poly-secondary}, we do not have any information about the relation between the secondary polytopes of~$\cA$ and $P$ as constructed in the proof of Proposition~\ref{prop:pc-to-poly}. So, given a point configuration~$\cA$, by using one of our two results we can either get a polytope with the same secondary polytope as~$\cA$ or a polytope with a tight span isomorphic to one of the tight spans of~$\cA$ but in general not both.
\begin{rem} 
\begin{enumerate}
\item Proposition~\ref{prop:pc-to-poly} enables us to give examples of $d$"=dimensional point configurations with tight spans equal to tight spans of $(d+1)$"=dimensional polytopes. Especially, examples of coarsest subdivisions of a point configuration whose tight spans have a given property directly give examples of coarsest subdivisions of a polytopes whose tight spans have the same property. We will make heavy use of this in the sequel, especially because this allows us to have the examples in lower dimension.
\item  Although coarsest subdivisions are mapped to coarsest subdivisions via the construction in the proof of Proposition~\ref{prop:pc-to-poly}, starting with a regular triangulation of a point configuration, we normally do not arrive at a triangulation of the polytope. For example, consider the point configuration~$\cA$ from Example~\ref{ex:square} and the lifting function $w=(-1/2,-1/2,-1/2,-1/2,-1)$. The subdivision $\subdivision w \cA$ is the triangulation depicted in the right part of Figure~\ref{fig:square-pc}. The polytope $P$ constructed in the proof of Proposition~\ref{prop:pc-to-poly} has ten vertices and the subdivision~$\subdivision {w'} P$ has four maximal cells that have six vertices and are combinatorially isomorphic to prisms over simplices.
\end{enumerate}
\end{rem}

\subsection{Existence of Tight Spans with the G"=property}

When considering tight spans, one might wonder which polytopal complexes might arise as the tight span of some regular subdivision of a polytope (or a point configuration). We will now give an answer for this question in the special case where the subdivision has the G"=property: In this case, where the tight span is a single polytope, it can be any polytope.

\begin{thm}\label{thm:ts-all}
Let $P$ be a $d$"=dimensional polytope with $n$ vertices. Then there exists a $d+1$"=dimensional polytope $P'$ with $2(n+1)$ vertices and a regular subdivision $\subdivision w {P'}$ of $P'$  such that the tight span $\tightspan{w}{P'}$ is affinely isomorphic to $P$.
\end{thm}

For the proof we need some notions about polytope polarity. We only give the notions and results we use here and refer the reader to \cite[Section~2.3]{MR1311028} or \cite[Section~3.4]{MR1976856} for details.

For a set $A\subset \RR^{d}$, the \emph{polar set} $A^\circ$ is defined as
\[
A^\circ\ =\ \SetOf{y\in \RR^{d}}{\scp xy\leq 1}.
\]
If $A$ is a compact convex set (e.g., a polytope) with $0\in\intt A$ then $(A^\circ)^\circ=A$. For a  polytope $P$ with $0\in\intt P$ (Note that this implies that $P$ is $d$"=dimensional.), the polar~$P^\circ$ equals $\conv(\vertt P)^\circ$ and is also a $d$"=dimensional polytope with $0\in\intt P^\circ$, called the \emph{polar} (or \emph{dual}) of~$P$. The face lattices of $P$ and $P^\circ$ are anti"=isomorphic.

\begin{proof}[Proof of Theorem~\ref{thm:ts-all}]
We assume that $P\subset\RR^d$ is $d$"=dimensional and that $0\in\intt P$, and we denote by $v_1,\dots,v_n$ the vertices of $P^\circ\subset\RR^{d}$.

Define the point configuration~$\cA\subset\RR^{d}$ as~$\cA=\{-v_1,\dots,-v_n,0\}$, and the lifting function $w:\cA\to \RR$ by $w(-v_i)=1$ for all $i\in\{1,\dots,n\}$, and $w(0)=0$. (Since $0$ is in the interior of $\conv\{-v_i\}\isom P^\circ$ the subdivision $\subdivision{w}{\cA}$ is obtained by coning from~$0$.) We get that
\begin{align*}
\envelope w \cA\ &=\ \SetOf{x\in\RR^{d+1}}{\left(\begin{smallmatrix} 1 & -v_1\\\vdots&\vdots\\1&-v_n\\1&0\end{smallmatrix}\right) x\geq -\left(\begin{smallmatrix} 1\\\vdots\\1\\0\end{smallmatrix}\right)}\\
&=\ \SetOf{(x_1,x')\in\RR_{\geq 0}\times\RR^d}{-x_1+\scp{v_i}{x'}\leq 1 \text{ for all } 1\leq i\leq \card \cA}\\
&=\ \RR_{\geq 0}\times (P^\circ)^\circ\,.
\end{align*}
This implies that $\tightspan w \cA=\{0\}\times(P^\circ)^\circ=\{0\}\times P$. So we have constructed a point configuration $\cA\subset \RR^d$ with $n+1$ points and a regular subdivision $\subdivision w \cA$ of $\cA$ such that $\tightspan w \cA$ is isomorphic to $P$. By Proposition~\ref{prop:pc-to-poly}, this implies the existence of a $d+1$"=dimensional polytope $P'$ with $2(n+1)$ vertices and a regular subdivision $\subdivision {w'} {P'}$ such that $\tightspan {w'} {P'}$ is isomorphic to $P$.
\end{proof}

\section{$k$"=Splits}\label{sec:k-splits}

We will now start our investigation of the coarsest subdivisions of a point configuration~$\cA$. The motivation of our definition is the notation of \emph{split} of a polytope defined in~\cite{MR2502496}. A split is a coarsest subdivision with exactly two maximal faces. It has the property that it contains exactly one interior face of  codimension one. This is the starting point of our generalization. We call a coarsest subdivision $\subdiv$ of~$\cA$ with $k$~maximal faces a \emph{$k$"=split} if $\subdiv$ has an interior face of codimension $k-1$.

It is easily seen that $\subdiv$ is a $k$"=split if and only if the tight span $\tightspan{\subdiv}{\cA}$ is a $(k-1)$"=dimensional simplex. So, in particular, all $k$"=splits have the G"=property.

\subsection{$1$"=Splits}\label{ssec:1splits}

For polytopes, 2"=splits are the \qm{simplest} possible non"=trivial subdivisions. However, general point configurations can have even simpler subdivisions: the $1$"=splits. For example, in the point configuration of Example~\ref{ex:square} (see Figure~\ref{fig:square-pc}), the subdivision with the sole maximal cell $\{1,2,3,4\}$ is non"=trivial. In general, for any $p\in\conv(\cA\setminus \{p\})$ there exists a subdivision $S^p$ of~$\cA$ with the unique maximal face~$\cA \setminus \{p\}$. (This includes configurations in convex position where one of the points occurs several times.) So for a point configuration~$\cA$ there does not exist a $1$"=split if and only if there exists a polytope $P$ such that~$\cA=\cA(P)$.

\begin{rem}\label{rem:oriented-matroid}
By the definition of 2"=split of a point configuration, it is clear that the set of 2"=splits of a point configuration~$\cA$ only depends on the oriented matroid of~$\cA$ as for polytopes; see~\cite[Remark~3.2]{MR2502496}. This is also obviously true for $1$"=splits.
\end{rem}

Given a $1$"=split $S^p$ of~$\cA$, we define a lifting function $w_p$ by $w_p(p)=1$ and $w_p(a)=0$ for all $a\in\cA$ with $a\not=p$. This obviously induces $S^p$. So all $1$"=splits are regular subdivisions. It is easily seen, that the tight span $\tightspan{w_p}{\cA}$ of any $1$"=split $S^p$ only consists of the single point $(0,\dots,0)$.

\subsection{Splits and the Split Decomposition}\label{sec:2splits}

A split of a polytope $P$ is a decomposition $S$ of~$P$ with exactly two maximal cells. So the splits of $P$ are the $2$"=splits of the point configuration~$\cA(P)$. Similarly, for a point configuration~$\cA$, we will define a \emph{split} of~$\cA$ as a $2$"=split of~$\cA$. Note that in the definition of split of a polytope it is not necessary to require that $S$ is a coarsest subdivision. However, the following example shows that this is needed for point configurations.

\begin{exmp}\label{ex:square2}
Let~$\cA$ be the point configuration from Example~\ref{ex:square}; see Figure~\ref{fig:square-pc}. Consider the subdivision $\subdiv_1$ with maximal cells $\{1,2,3,5\}$ and $\{1,3,4,5\}$ and the subdivision $\subdiv_2$ with maximal cells $\{1,2,3\}$ and $\{1,3,4\}$. Only $\subdiv_1$ is a 2"=split of~$\cA$ since $\subdiv_2$ is coarsened by the $1$"=split $S^5$.
\end{exmp}

The reason for this difference is that point configurations may have $1$"=splits, whereas polytopes may not. However, we have the following characterization of $2$"=splits of point configurations, whose simple proof we omit.

\begin{lem}
Let $S$ be a subdivision of~$\cA$ with exactly two maximal faces~$S_+$ and~$S_-$. Then the following statements are equivalent.
\begin{enumerate}
\item $S$ is a 2"=split of~$\cA$,
\item $S$ is a coarsest subdivision of~$\cA$,
\item $S_+=\conv S_+ \cap \cA$ and $S_-=\conv S_- \cap \cA$.
\end{enumerate}
\end{lem}

For a 2"=split~$S$ of a polytope $P$, there exists a hyperplane $H_S$ that defines $S$, and a hyperplane $H$ (that meets the relative interior of $P$) defines a 2"=split if and only if it does not meet any edge of $P$ in its relative interior. As well, for a 2"=split $S$ of a point configuration~$\cA$, there exists a hyperplane~$H$ inducing a 2"=split. However, the condition has to be modified a bit: A hyperplane $H$ defines a 2"=split of~$\cA$ if and only if it meets $\conv \cA$ in its interior and for all edges $E$ of $\cA$ we have that $H\cut E$ is either empty, a point of~$\cA$, or~$E$ itself. Here an \emph{edge}~$E$ of a point configuration $\cA$ is defined as the convex hull of two points in $\cA$ that are contained in some edge of the polytope $\conv \cA$. %We give a complete generalization of Observation~\ref{obs:facets} and Proposition~\ref{prop:split-condition} to point configurations. One can give similar generalizations of the results of Section~\ref{sec:split-conds}, for example for the conditions of compatibility, and so on.
This leads to the following statement 
% \begin{prop}\label{prop:split-condition-pc}
% Let~$\cA$ be a point configuration and $H$ a hyperplane that does not intersect $\conv \cA$ in the interior. Then the following are equivalent.
% \begin{enumerate}
% \item  $H$ induces a 2"=split on~$\cA$,
% \item $H$ meets all edges $E$ of~$\cA$ in an element of~$\cA$, $E$, or $\emptyset$,
% \item $H$ meets all faces of~$\cA$ in a face of~$\cA$ or induces a 2"=split on them,
% \item $H$ meets all facets of~$\cA$ in a face of~$\cA$ or induces a 2"=split on them,
% \item all vertices of the subdivision of~$\cA$ with maximal faces~$\cA\cut H_+$ and~$\cA\cut H_-$ are elements of~$\cA$,
% \item $H\cut \conv \cA=\conv(\cA \cut H)$.
% \end{enumerate}
% \end{prop}
% Again, $H_+$ and $H_-$ denote the two halfspaces in which $\RR^{d+1}$ is divided by the hyperplane~$H$.
%
% An immediate consequence of this characterization is the following lemma 
which says that by adding points in the convex hull one cannot lose 2"=splits.

\begin{lem}\label{lem:lose-splits}
Let~$\cA$,~$\cA'$ be point configurations with $\cA'\subset \cA$ and $\conv \cA=\conv \cA'$. If~$S$ is a 2"=split of ~$\cA'$ with maximal faces $S_+$ and $S_-$, then~$\cA$ has a 2"=split $S'$ with maximal faces $S'_+=\conv S_+\cap \cA$ and $S'_-=\conv S_-\cap \cA$.

Especially, if $S$ is a 2"=split of the polytope $\conv \cA$ with maximal faces $S_+$ and $S_-$, then~$\cA$ has a 2"=split $S'$ with maximal faces $S'_+= S_+\cap \cA$ and $S'_-= S_-\cap \cA$.
\end{lem}

\begin{rem} A point configuration is called \emph{unsplittable} if it does not admit any 2"=split. It follows from
Lemma~\ref{lem:lose-splits} that a two"=dimensional point configuration~$\cA$ with $\conv \cA$ not being a simplex
cannot be unsplittable. But -- in contrast to the polytope case -- there are a lot of different point configurations
whose convex hulls are simplices. In fact, such a point configuration~$\cA$ is unsplittable if and only if there is no
$a\in\cA$ which is in the relative interior of an edge of $\conv \cA$. This gives us a lot of non"=trivial unsplittable
two"=dimensional point configurations, namely all point configurations having a point in the relative interior but no point in the
relative interior of an edge. So the simplest non"=trivial unsplittable point configuration is a triangle with a point
in its interior. \end{rem}

For point configurations, we have the following generalization of the Split Decomposition Theorem \cite[Theorem 2]{MR1153934}, \cite[Theorem~3.10]{MR2502496}, \cite[Theorem~2.2]{MR2252108}. A lifting function $w:\cA\to \RR$ is called \emph{split prime} if the subdivision $\subdivision w \cA$ is not refined by any $1$"=split or $2$"=split.

\begin{thm}[Split Decomposition Theorem for Point Configurations]\label{thm:splitdecomposition-pc}
 
  Let $\cA$ be a point configuration. Each weight function $w:\cA\to \RR$ has a coherent decomposition
\begin{equation}\label{eq:splitdecomposition-pc} w \ = \ w_0 + \sum_{\text{$S^p$ a $1$"=split of~$\cA$}} \lambda_{p} w_p
+\sum_{\text{$S$ 2"=split of~$\cA$}} \lambda_S w_S\, , \end{equation} where $w_0$ is 2"=split prime, and this is unique
among all coherent decompositions of~$w$ into $1$"=splits, $2$"=splits, and a split weight function.
\end{thm}
 \begin{proof} The proof works in the same manner as the proof of
\cite[Theorem~3.10]{MR2502496}. We first consider the special case where the subdivision $\subdivision w \cA$ is a
common refinement of $1$"=splits and $2$"=splits. The $1$"=splits coarsening $\subdivision w \cA$ are those $S^p$ where
$p$ is not contained in any face of $\Asubdivision w\cA$. Moreover, each face $F$ of codimension~$1$ in
$\subdivision{w}{\cA}$ defines a unique split $S$ whose split hyperplane is $\aff F$. Whenever $S$ is an arbitrary
split of~$\cA$, then there exists some $\lambda_S>0$ such that $(w-\lambda_S w_S)+\lambda_S w_S$ is coherent if and only if $H_S\cap \cA$ is a face of $\subdivision w \cA$ of codimension
one.  So we get a coherent decomposition $w =\sum_{p\in \cA} \lambda_p w_p +\sum_S \lambda_S w_S$, where the
second sum ranges over all splits $S$ of $\cA$.  Note that the uniqueness follows  from the fact that for each
codimension"=one"=face of $\subdivision w \cA$ there is a unique split~$S$ whose split hyperplane~$H_S$ contains it.

For the general case, we define
\[
w_0 \ := \  w - \sum_{\text{$S$ split of $\cA$}} \lambda_S w_S-\sum_{\text{$S^p$ a $1$"=split of~$\cA$}} \lambda_p w_p \, .
  \]
  This weight function is split prime by construction, and the uniqueness of the split decomposition of $w$
  follows from the uniqueness of the split decomposition of $w-w_0$.
\end{proof}

\subsection{General $k$"=Splits}

\begin{exmp}
An example of a $k$"=split is given by taking a $(k-1)$"=dimensional simplex with a point in the interior and coning from that point. For an example of a polytope (with less vertices than that one could obtain from Proposition~\ref{prop:pc-to-poly}), one can take a bipyramid over a $(k-1)$"=dimensional simplex and cone from the edge connecting the two pyramid vertices.
\end{exmp}

\begin{figure}[tb]
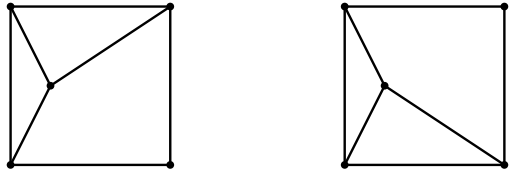
\centering
  \includegraphics[scale=1]{3-splits.0}\hspace{2cm}
  \includegraphics[scale=1]{3-splits.1}
  \caption{A point configuration with a codimension"=two"=face (the interior point) that corresponds to two different $3$"=splits.}
  \label{fig:3-splits}
\end{figure}

We know that to each 2"=split $S$ there corresponds a unique hyperplane $H_S$ that defines~$S$. For general $k$"=splits, it is still true that to a $k$"=split $\Sigma$ (for $k>1$) there corresponds a unique subspace of codimension $k-1$. However, for 2"=splits we also have the property that if a hyperplane $H$ defines a 2"=split, this 2"=split is uniquely determined by $H$. This does not hold any more for $k$"=splits with $k\geq 3$; see Figure~\ref{fig:3-splits}.% and Example~\ref{ex:3cube:3}.

In Section~\ref{sec:2splits}, we have seen that a hyperplane $H$ defines a 2"=split of a point configuration if and only if it meets all edges $E$ of~$\cA$ in an element of~$\cA$, $E$, or the empty set. One direction of this generalizes to $k$"=splits as follows.

\begin{prop}\label{prop:k-split-condition}
If $U$ is the unique codimension"=$(k-1)$"=subspace of $\aff \cA$ corresponding to some $k$"=split of a point configuration~$\cA$, then the following equivalent conditions are satisfied.
\begin{enumerate}
\item $U$ meets all faces $F$ of~$\cA$ with $\dim F\leq k-1$ in a face of~$\cA$ or corresponds to an $l$"=split of them with $l\leq k$,\label{cond:k-split-cond-1}
\item $U$ meets all faces of~$\cA$ in a face of~$\cA$ or corresponds to an $l$"=split of them for some  $l\leq k$,\label{cond:k-split-cond-2}
\item $U$ meets all facets of~$\cA$ in a face of~$\cA$ or corresponds to an $l$"=split of them for some $l\leq k$.\label{cond:k-split-cond-3}
\end{enumerate}
\end{prop}
\begin{proof}
First one sees that if $\subdiv$ is a $k$"=split of~$\cA$, the induced subdivision to each face of~$\cA$ has to be an $l$"=split for some $l\leq k$ or the trivial subdivision. This implies that all conditions have to be satisfied. That \eqref{cond:k-split-cond-1} implies \eqref{cond:k-split-cond-2} follows from the fact that if a codimension"=$(k-1)$"=subspace $U$ intersects some face $F$ with $\dim F\geq k$ in its interior, the subspace $U$ has to intersect some of the faces of $F$ of dimension $k-1$. That \eqref{cond:k-split-cond-3} is also equivalent follows by applying the equivalence of \eqref{cond:k-split-cond-1} and \eqref{cond:k-split-cond-2} to~$\cA$ and its facets.
\end{proof}

\begin{figure}[htb]\centering
  \includegraphics[scale=0.45]{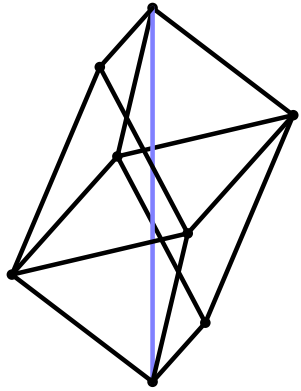}
  \caption{A polytope with an interior edge that does not correspond to a $3$"=split.}
  \label{fig:no-3-split}
\end{figure}

However, in contrast to the 2"=split case, the converse of Proposition~\ref{prop:k-split-condition} does not hold if $k\geq 3$. For an example, consider the polytope depicted in Figure~\ref{fig:no-3-split}. The codimension"=two"=subspace spanned by the top and bottom vertices does not correspond to any $3$"=split.

A key property of 2"=splits \cite[Lemma~3.5]{MR2502496} is shared by $k$"=splits: They are regular subdivisions.

\begin{figure}[htb]
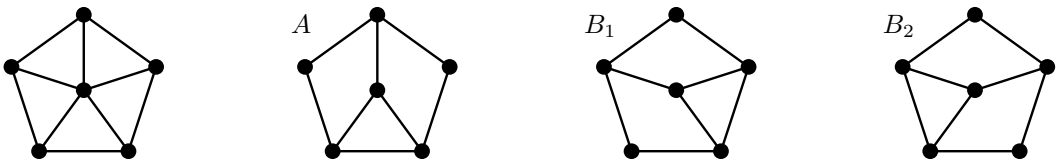
\centering
  \includegraphics[scale=1]{no-sdt.0}\hspace{1.5cm}
  \includegraphics[scale=1]{no-sdt.1}\hspace{1.5cm}
  \includegraphics[scale=1]{no-sdt.2}\hspace{1.5cm}
  \includegraphics[scale=1]{no-sdt.3}
  \caption{There is no unique $3$"=split decomposition.}
  \label{fig:no-sdt}
\end{figure}

\begin{thm}\label{thm:k-regular}
 All $k$"=splits are regular.
\end{thm}
\begin{proof}
  Let~$\cA\subset\RR^d$ be a $d$"=dimensional point configuration and $\subdiv$ a $k$"=split of~$\cA$. Then $\subdiv$ has a unique interior face $F$ such that $U:=\aff F$ has dimension $d-(k-1)$. We can assume without loss of generality  that the origin is contained in $\conv F$. Let now $\pi$ be the projection orthogonal to $U$. We consider the subdivision $\subdiv':=\pi(\subdiv)$ of the $(k-1)$"=dimensional point configuration~$\cA':=\pi(\cA)$ with the origin as an interior vertex. If we now take for each face $F$ of $\subdiv'$ the cone spanned by $F$, we get a polyhedral fan $\cF$ subdividing $\RR^{k-1}$. The dual complex of $\cF$ is isomorphic to $\tightspan {\subdiv'}{\cA'}$ and hence to $\tightspan \subdiv \cA$. For each of the $k$ rays $r_i$ of this fan (which correspond to interior faces of dimension $d-k+2$ of $\subdiv$), we take a vector $v_i$ of length one that spans this ray. Each point $a'\in\cA'$ is contained in the relative interior of a unique cone $C\in\cF$ and can uniquely be written as $a'=\sum_{i=1}^{k}
 \lambda^{a'}_i v_i$ where $\lambda^{a'}_i\geq 0$ for all $1\leq i\leq k$ and $\lambda^{a'}_i>0$ if and only if $v_i\in C$. Now we define a weight function $w_\subdiv:\cA\to \RR$ via $w_\subdiv(a):=\sum_{i=1}^{k-1}\lambda^{\pi(a)}_i$. This lifting function $w_\subdiv$ defines $\subdiv$.
\end{proof}
\begin{rem}
One might ask whether there exists some generalization of the Split Decomposition Theorem \ref{thm:splitdecomposition-pc} to $k$"=splits. However, even if one fixes some $k\geq 3$, no similar result can be valid: The triangulation to the left of Figure~\ref{fig:no-sdt} can be obtained as the common refinement of the $3$"=split $A$ and either of the two $3$"=splits~$B_1$ and~$B_2$.
\end{rem}

\subsection{Approximation of Secondary Polytopes}

As explained in Section~\ref{sec:secondary-polytope}, the facets of the secondary polytope of a point configuration $\cA$ are in bijection with the coarsest regular subdivisions of~$\cA$ and the facet"=defining inequalities can be explicitly computed from weight functions for that subdivisions. Hence each $k$"=split of $\cA$ gives rise to such an inequality.

Firstly, we are only interested in the $1$- and $2$"=splits, for which weight functions are computed very easily. As in the case of a polytope, we can define the \emph{split polyhedron} $\SplitPoly{\cA}$ of a point configuration~$\cA$. It a $(\card \cA-d-1)$"=dimensional polyhedron in $\RR^{\card\cA}$ defined by one inequality for each $1$- or $2$"=split together with a set of equations defining the affine hull of $\SecondaryPolytope \cA$. Remark~\ref{rem:oriented-matroid} shows that the split polyhedron only depends on the oriented matroid of~$\cA$ and hence can be seen as a common approximation of the secondary polytope of all point configurations with the same oriented matroid.

The proof  of Theorem~\ref{thm:k-regular} allows us to generalize this to arbitrary $k$"=splits: For each $k$"=split $\Sigma$ of~$\cA$ we construct a weight function $w_\Sigma$ as in the proof of Theorem~\ref{thm:k-regular}.  We get an explicit description of the inequality~$I_\subdiv$ defining the corresponding facet of $\SecondaryPolytope \cA$. The \emph{$k$"=split polyhedron} $k-\SplitPoly \cA$ of~$\cA$ is then defined as the intersection of  $\SplitPoly \cA$ with all halfspaces defined by some $I_\Sigma$ where $\Sigma$ ranges over all $l$"=splits of~$\cA$ with $3\leq l\leq k$.

This gives us a descending sequence of outer approximations for $\SecondaryPolytope \cA$. Obviously, since a $d$"=dimensional point configuration cannot have any $k$"=splits for $k>d+1$, this sequence eventually becomes constant at the value $(d+1)-\SplitPoly \cA$. If $\cA=\cA(P)$ for some polytope $P$, then $P$ cannot have $k$"=splits for $k>d$, so the sequence becomes already constant at the value $d-\SplitPoly P$. This is the best possible approximation of the secondary polytope that one may obtain via $k$"=splits.

\subsection{Totally $k$"=Splittable Point Configurations}\label{sec:tkspc}

In \cite{MR2502496}, a polytope~$P$ was defined to be \emph{totally splittable} if and only if all regular subdivisions of~$P$ are refinements of splits or, equivalently, if and only if $\SecondaryPolytope P=\SplitPoly P$. These polytopes can be completely classified \cite[Theorem~9]{MR2502496}: A polytope $P$ is totally splittable if and only if it has the same oriented matroid as a simplex, a cross polytope, a polygon, a prism over a simplex, or a (possibly multiple) join of these polytopes. We generalize this definition from polytopes to point configurations and from $2$"=splits to $k$"=splits for arbitrary $k$. A point configuration~$\cA$ is called \emph{totally $k$"=splittable} if and only if $\SecondaryPolytope \cA=k-\SplitPoly \cA$. This is equivalent to saying that all regular subdivisions of~$\cA$ are common refinements of $l$"=splits with $l\leq k$.

So the totally $k$"=splittable point configurations are those point configuration whose secondary polytopes can be entirely computed by computing the $l$"=splits for all $l\leq k$ and then constructing the weight functions as in the proof of Theorem~\ref{thm:k-regular}.

 Before closing this section with some examples of totally $k$"=splittable point configurations, we remark that totally $k$"=splittable polytopes obviously have the G"=property, since all $k$"=splits have the G"=property.

\begin{exmp}\label{ex:3cube:3}
The $3$"=cube $C_3$ has 14 2"=splits (see \cite[Example~3.8]{MR2502496}), and eight $3$"=splits: Each diagonal of the cube corresponds to two $3$"=splits by subdividing $C_3$ into three square pyramids with one of the vertices of the diagonal as apex. In particular, $C_3$ is not totally $2$"=splittable. By using the 14 inequalities obtained for the weight functions defining the $2$"=splits and the eight inequalities obtained from the weight functions for defining eight $3$"=splits, we can compute $3-\SplitPoly {C_3}$. It is easily observed that all triangulations of $C_3$ are obtained as refinements of $2$"=splits or $3$"=splits, so $C_3$ is totally $3$"=splittable. This gives us a new computation of the secondary polytope of the $3$"=cube, verifying the results of Pfeifle~\cite{Pfeifle00}.
\end{exmp}

\begin{exmp}
 The secondary polytope of the $4$"=cube $C_4$ was computed by Huggins, Sturmfels, Yu, and Yuster~\cite{Grier+08}. It has $80,876$ facets that come in $334$ orbits. An inspection of their results shows that four of these orbits are 2"=splits, five are $3$"=splits, and three are $4$"=splits. So $C_4$ is not totally $4$"=splittable, hence not totally $k$"=splittable for any $k$.
\end{exmp}

\begin{prop}
Let $\cA$ be a $(\card \cA -2)$"=dimensional point configuration. Then $\cA$ is totally $(\card \cA-1)$"=splittable. If $\cA$ is the vertex set of a polytope, then $\cA$ is totally $(\card\cA-2)$"=splittable.
\end{prop}
\begin{proof}
Since $\cA$ is $(\card \cA -2)$"=dimensional, the Gale dual $\cB$ of $\cA$ is one"=dimensional. So the maximal faces of $\ChamberComplex \cB$ are the two rays $\pos 1$ and $\pos (-1)$. By Theorem~\ref{thm:secondaryDuality}, the sole non"=trivial subdivisions of $\cA$ are a $k$"=split and an $l$"=split, where $k$ is the number of $b\in\cB$ with $\pos \cB = \pos 1$ and $l$ is the number of $b\in \cB$ with $\pos b=\pos (-1)$. Since $\cB$ positively spans the whole space, we have $k,l\geq 1$. If $\cA$ is the vertex set of a polytope, we have $k,l\geq 2$ by \cite[Section~5.4, Theorem~2]{MR1976856}. The fact that $k+l\leq \card \cA$ then shows the claim.
\end{proof}

\section{General Coarsest Subdivisions}\label{sec:k-subdivs}

Now we will discuss coarsest subdivisions of point configurations that are not necessarily $k$"=splits. To simplify the notation, we call a coarsest subdivision with $k$ maximal faces a \emph{$k$"=subdivision}.

For $1$"=subdivisions and $2$"=subdivisions, it is easily seen that their tight spans are points and line segments, respectively. Especially, all $1$"=subdivisions are $1$"=splits and all $2$"=subdivisions are 2"=splits. We will see in Lemma~\ref{lem:3-subdivisions} that $3$"=subdivisions are $3$"=splits, too. However, for $k$"=subdivisions with $k> 3$ the tight spans get much more complicated. We will investigate these tight spans in this section. First, we give two general statements about the tight spans of $k$"=subdivisions. Note that everything we prove in this section is not only true for regular subdivisions but also for non"=regular subdivisions and their tight spans as defined in Remark~\ref{rem:non-regular}.

By Theorem~\ref{thm:ts-all}, for each polytope $P$ there exists some polytope $P'$ whose tight span is (the complex of faces of)~$P$. The next proposition shows that this is not true if one only considers $k$"=subdivisions, that is, coarsest subdivisions.

\begin{prop}\label{prop:k-no-n-gon}
 Let~$\cA$ be a point configuration, $k>3$, and $\subdiv$ a $k$"=subdivision of~$\cA$. Then the tight span $\tightspan{\subdiv}{\cA}$ is not a $k$"=gon.
\end{prop}

\begin{proof}
Suppose we have some  subdivision $\subdiv$ of~$\cA$ whose tight span is a $k$"=gon. The $k$"=gon corresponds to some codimension"=two"=face $F$ of $\subdiv$. The facets of $F$ are all contained in the boundary of $\conv \cA$ since any facet of $F$ that is an interior face would correspond to a three"=dimensional face of $\tightspan{\subdiv}{\cA}$. So we have $F=\aff F\cut \cA$. The edges of the $k$"=gon are dual to codimension"=one"=faces of $\subdiv$ whose intersection is $F$. Call these faces $F_1,\dots,F_k$ (We consider the indices modulo $k$.), where $F_1$ is chosen arbitrary and the others are numbered in counter"=clockwise order. Furthermore, the maximal cell of~$\subdiv$ between $F_i$ and $F_{i+1}$ is called $C_i$. For each cell $C_i$ one can measure the angle $\alpha_i$ between the (hulls of the) two consecutive faces $F_i$ and $F_{i+1}$. Obviously, $\sum_{i=1}^k \alpha_i = 2\pi$, and, since $k>3$, there exists at least one $i$ with $\alpha_i+\alpha_{i+1}\leq \pi$.

We now distinguish two cases. If $\alpha_i+\alpha_{i+1} = \pi$, the hyperplane $\aff F_i=\aff F_{i+2}$ defines a 2"=split of~$\cA$ refined by $\subdiv$, contradicting the fact that $\subdiv$ was supposed to be a coarsest subdivision. On the other hand, $\alpha_i+\alpha_{i+1} < \pi$ implies that $\conv C_i \union \conv C_{i+1}$ is convex. Therefore, we can construct a new subdivision $\subdiv'$ of~$\cA$ with the $k-1$ maximal faces $C_1,\dots,C_{i-1},C_i\nobreak\union\nobreak C_{i+1},C_{i+2},\dots,C_k$. Since $\alpha_i+\alpha_{i+1} < \pi$, the faces $F_i,F_{i+2}\in\subdiv$ are also faces of $\subdiv'$, what ensures that (SD3) holds, and hence $\subdiv'$ is a valid subdivision of~$\cA$.
\end{proof}

Note that this only shows that $k$"=gons with $k>3$ cannot be the sole maximal cell of the tight span. It can well be that a polygon occurs as a maximal cell of a tight span of a $k$"=subdivision if there are other maximal cells. For the simplest example see the top left part of Figure~\ref{fig:5-subdiv-2}.

\begin{prop}\label{prop:1-connected}
Let~$\cA$ be a point configuration and $\subdiv$ a $k$"=subdivision of~$\cA$. Then the graph of the tight span $\tightspan{\subdiv}{\cA}$ is $2$"=connected, that is, it is still connected if one removes any vertex.
\end{prop}
\begin{proof}
We will show that for a subdivision $\subdiv$ of~$\cA$ for which the graph of its tight span is not $2$"=connected there exists a subdivision $\subdiv'$ of~$\cA$ that coarsens $\subdiv$.

So suppose that there exists a vertex $v$ of $\tightspan{\subdiv}{\cA}$ such that $\tightspan{\subdiv}{\cA}\setminus\{v\}$ is not connected. Let $T$ be the set of vertices of some connected component of  $\tightspan{\subdiv}{\cA}\setminus\{v\}$. For a vertex~$w$ of $\tightspan{\subdiv}{\cA}$ the corresponding maximal cell of $\subdiv$ is denoted by $w^\circ$ . We then define a new subdivision~$\subdiv'$ of~$\cA$ by deleting all maximal cells $w^\circ$ with $w\in T\union\{v\}$ and adding $F:=\Union_{w\in T\union\{v\}}w^\circ$ as a new maximal cell of~$\subdiv'$. In order to show that $\subdiv'$ is actually a subdivision of~$\cA$, we have to show that (SD3) holds. 

We first show that $C:=\Union_{w\in T\union\{v\}}\conv w^\circ$ is convex. So assume that there exists $x,y\in \relint C$ such that the line segment $l$ connecting $x$ and $y$ is not entirely contained in $C$. Then $l$ has to intersect two codimension"=one"=cells $C_1$ and $C_2$ of those remaining in $\subdivG$. However, by our assumption that $T$ is the set of vertices of some connected component of $\tightspan{\subdiv}{\cA}\setminus\{v\}$, the edges of $\tightspan \subdiv \cA$ corresponding to these cells can only be connected to $v$. So $C_1$ and $C_2$ are facets of $\conv v^\circ$ and this implies that $\conv v^\circ$ is not convex, a contradiction.

To finish the proof of (SD3), note that an improper intersection cannot happen in the interior of $\conv \cA$ since all interior faces of $C$ are interior faces of $v^\circ$ by assumption. However, any improper intersection of faces $F_1,F_2$ in the boundary of $\conv \cA$ would yield an improper intersection of some interior faces $F_1', F_2'$ with $F_1\subset F_1',F_2\subset F_2'$. So~$\subdiv'$ is a subdivision of~$\cA$ that coarsens $\subdiv$, as desired.
\end{proof}
As a third condition for the tight span of a $k$"=subdivision, we note that any tight span of a regular subdivision has to be a contractible \cite[Lemma~4.5]{MR2233884} and hence simply"=connected polyhedral complex. It can be shown that this is true also for non"=regular subdivisions. Additionally, this leads to the following important corollary.

\begin{cor}\label{cor:k-no-edge}
 Let~$\cA$ be a point configuration and $\subdiv$ a coarsest subdivision of~$\cA$ that is not a 2"=split. Then all maximal faces of the polyhedral complex $\tightspan{\subdiv}{\cA}$ are at least two"=dimensional.
\end{cor}
\begin{proof}
Suppose there exists some edge $E$ in $\tightspan \subdiv \cA$ connecting $v$ and $w$ that is a maximal face. Since $\subdiv$ is not a 2"=split, we can assume that one of the vertices of $E$ is strictly contained in another face of $\tightspan \subdiv \cA$. If we delete this vertex from $\tightspan \subdiv \cA$, by Proposition~\ref{prop:1-connected}, the remainder is still connected. However, this implies that there has to be a path in the graph of $\tightspan \subdiv \cA$ connecting $v$ with $w$ without using $E$. This contradicts the simple connectedness.
  % Suppose that there exists some edge $E$ in $\tightspan{\Delta}{\cA}$ which is not contained in some higher dimensional face. The edge $E$ corresponds to some codimension-1"=face $F$ of $\Delta$. The facets of $F$ are all contained in the boundary of $\conv \cA$ since any other facet of $F$ would correspond to a $2$"=dimensional face $E'$ of $\tightspan{\Delta}{\cA}$ with $E\subset E'$. So we have $\conv F=\aff F \cut \conv \cA$. Furthermore, $F$ is adjacent to exactly to full"=dimensional cells $F_1$, $F_2
%   So we have constructed a refinement $\Delta'$ of $\Delta$ which is non"=trivial since $k\geq 3$. This shows that $\Delta$ cannot be a $k$"=subdivision.
\end{proof}

\subsection{Tight Spans of $k$"=subdivisions for Small $k$}

Now, we will examine the tight spans of $k$"=subdivision for small $k$. We start out with a complete characterization of tight spans of $k$"=subdivisions for $k=3,4$.

\begin{lem}\label{lem:3-subdivisions}
Let~$\cA$ be a point configuration, and $\subdiv$ a $3$"=subdivision of~$\cA$. Then the tight span of $\subdiv$ is a triangle.
\end{lem}
\begin{proof}
 Obviously, the only simple connected polyhedral complexes with three points are a triangle or two line segments connected at one point. However, the latter cannot occur by Proposition \ref{cor:k-no-edge}.
\end{proof}

Since a $3$"=subdivision whose tight span is a triangle has an interior face of codimension $2$ we directly get.
\begin{cor}\label{cor:3-subdivs-3-splits}
All $3$"=subdivisions are $3$"=splits.
\end{cor}

\begin{figure}[htb]\centering
  \includegraphics[scale=0.6]{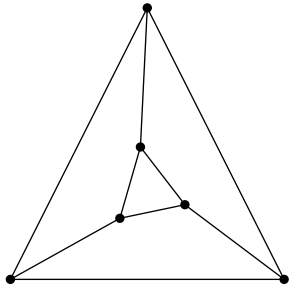}
  \caption{A non-regular $4$"=subdivision.}
  \label{fig:non-regular}
\end{figure}

\begin{rem}
Corollary~\ref{cor:3-subdivs-3-splits} and Theorem~\ref{thm:k-regular} imply that all $3$"=subdivisions and furthermore all subdivisions with at most three maximal faces are regular. This is not true anymore for subdivisions with four or more maximal faces. An example is the subdivision depicted in Figure~\ref{fig:non-regular}: Suppose that subdivision would be induced by a lifting function. One can assume that the three interior points are lifted to $0$. It is easily seen that one cannot choose the weights of the vertices of the outer triangle in such a way that the depicted subdivision is induced since the inner triangle is slightly rotated. This example is related to the so"=called \qm{mother of all examples} (of a non"=regular triangulation); see \cite[Section~7.1]{Triangulations}.
\end{rem}

\begin{lem}\label{lem:4-subdivs}
Let~$\cA$ be a point configuration, and $\subdiv$ a $4$"=subdivision of~$\cA$. Then the tight span of $\subdiv$ is either a tetrahedron, or it consists of three triangles with a common vertex, or it consists of two triangles glued together at one edge.
\end{lem}
\begin{proof}
We have to look at simply connected polyhedral complexes with four vertices. By Corollary~\ref{cor:k-no-edge}, we have the additional condition that all maximal cells have to be at least two"=dimensional. So the candidates are a tetrahedron, two triangles glued together at one edge, three triangles with a common vertex, or a quadrangle. However, the quadrangle cannot occur by Proposition~\ref{prop:k-no-n-gon}.
\end{proof}

\begin{figure}[htb]
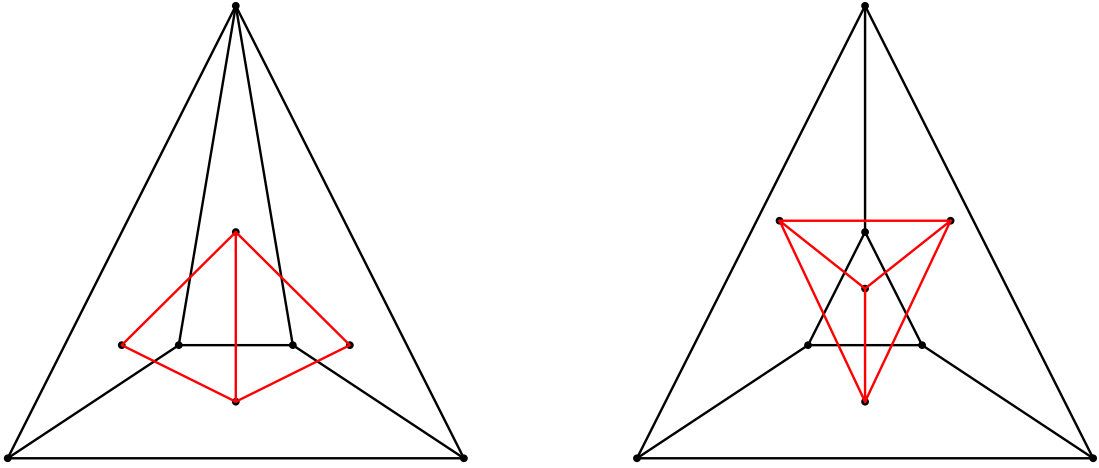
\centering
  \includegraphics[scale=1]{4-subdivs.0}\hspace{2cm}
  \includegraphics[scale=1]{4-subdivs.1}
  \caption{Two $4$"=subdivisions and their tight spans.}
  \label{fig:4-subdiv}
\end{figure}

\begin{exmp}\label{ex:4-subdiv-3d}
In Figure~\ref{fig:4-subdiv}, we depict examples of $4$"=subdivisions of point configurations together with their tight spans, which are the two two"=dimensional complexes from Lemma~\ref{lem:4-subdivs}. A $4$"=subdivision with a tetrahedron as tight span is a $4$"=split.
\end{exmp}

For $5$"=subdivisions, the number of possible tight spans gets much larger. However, we have here the first case of a simply connected polyhedral complex that cannot occur as a tight span of a $k$"=subdivision and is not excluded by Proposition~\ref{prop:k-no-n-gon} or Proposition~\ref{prop:1-connected}.
\begin{figure}[hp]
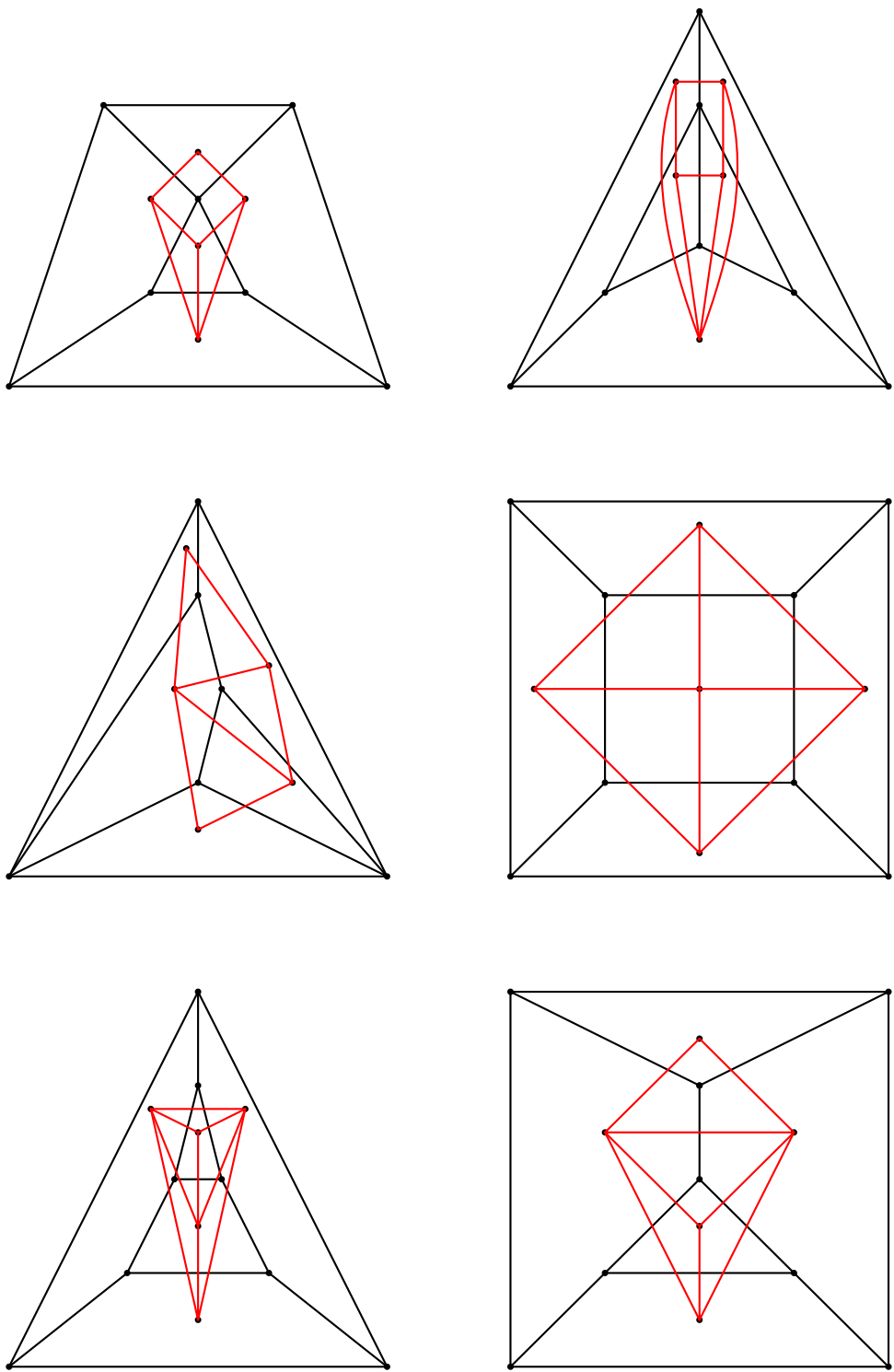
\centering
  \includegraphics[scale=0.9]{5-subdivs.1}\hspace{1.5cm}
  \includegraphics[scale=0.9]{5-subdivs.0}\vspace{1.5cm}
  \includegraphics[scale=0.9]{5-subdivs.2}\hspace{1.5cm}
  \includegraphics[scale=0.9]{5-subdivs.3}\vspace{1.5cm}
  \includegraphics[scale=0.9]{5-subdivs.4}\hspace{1.5cm}
  \includegraphics[scale=0.9]{5-subdivs.5}
  
  \caption{The $5$"=subdivisions with planar tight spans.}
  \label{fig:5-subdiv-2}
\end{figure}

\begin{lem}\label{lem:5-subdivs}
Let~$\cA$ be a point configuration and $\subdiv$ a $5$"=subdivision of~$\cA$. Then the tight span of $\subdiv$ cannot consist of a quadrangle and a triangle glued together at one edge.% is either a $4$"=dimensional simplex, a pyramid over a square, a bipyramid over a triangle, consists of two tetrahedra glued at a triangle, consists of a tetrahedron and a triangle glued at one edge, or one of the $2$"=dimensional complexes depicted in Figure~\ref{fig:5-subdiv-2}
\end{lem}
\begin{proof}
Suppose there exists a point configuration~$\cA$ and a subdivision $\subdiv$ of $\cA$ with such a tight span and let $E$ be the edge  of the tight span which is the intersection of the quadrangle and the triangle.We can now argue as in the proof of Proposition~\ref{prop:k-no-n-gon} by letting $F$ be the face of $\subdiv$ dual to the quadrangle. We adopt the notation from the proof of Proposition~\ref{prop:k-no-n-gon}. The only case that is not covered by the argument there is when the index $i$ is such that $C_i$ and $C_{i+1}$ are the cells corresponding to the vertices of~$E$ and $\alpha_i+\alpha_{i+1}< \pi$. However, in this case, we simply take $C:=C_i\union C_{i+1} \union C^\star$ instead of $C_i\union C_{i+1}$ as a new maximal cell, where $C^\star$ is the cell of $\subdiv$ corresponding to the unique non"=quadrangle vertex of the tight span. One now directly sees that (SD3) holds by the same argumentation as in the proof of Proposition~\ref{prop:k-no-n-gon}.
\end{proof}

\begin{exmp}
In Figure~\ref{fig:5-subdiv-2}, we depict examples of $5$"=subdivisions covering all planar tight spans that may occur. For the two topmost subdivisions it has to be carefully checked that these are really coarsest subdivisions, which is true because all unions of occurring cells are not convex.
\end{exmp}

\begin{figure}[hp]
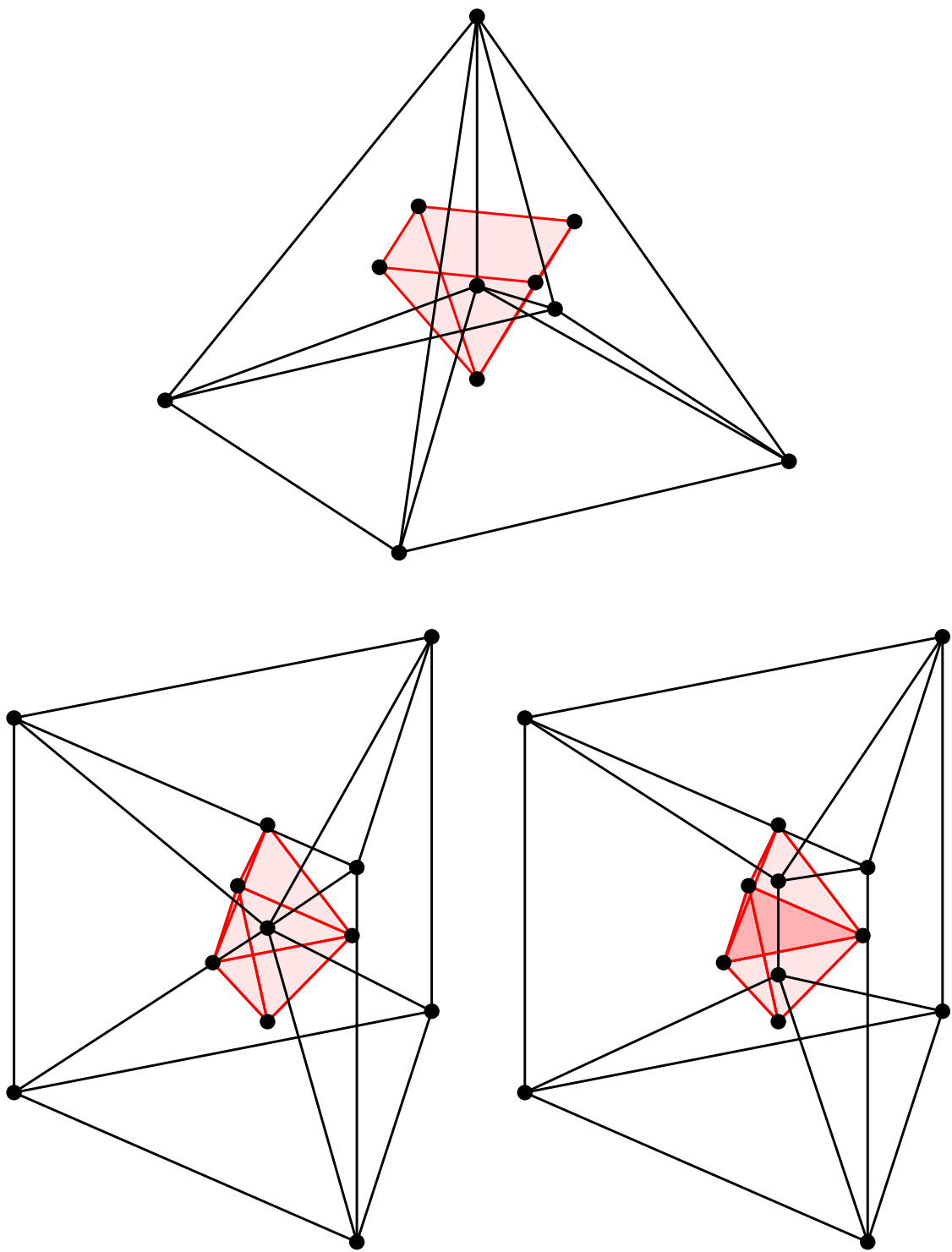
\centering
  \includegraphics[height=0.4\textheight]{5-subdivs-3d.0}\vspace{1cm}
  \includegraphics[height=0.45\textheight]{5-subdivs-3d.1}\hspace{1cm}
  \includegraphics[height=0.45\textheight]{5-subdivs-3d.2}
    \caption{Some $5$"=subdivisions with pure three"=dimensional tight spans.}
  \label{fig:5-subdiv-3}
\end{figure}

\begin{exmp}
In Figure~\ref{fig:5-subdiv-3}, we depict some examples of $5$"=subdivisions with pure three"=dimensional tight spans. The first tight span is a pyramid, and the subdivision is obtained by taking as point configuration the vertices of another pyramid $P$ together with any interior point $v$ and as maximal simplices the cones from $v$ over all facets of~$P$. (This is the same construction as in the proof of Theorem~\ref{thm:ts-all}; pyramids are self"=dual.) To the left, we have as tight span a bipyramid over a triangle, which is obtained in the same way by taking a prism over a triangle with one interior point.
The tight span of the subdivision to the right of Figure~\ref{fig:5-subdiv-3} consist of two tetrahedra glued at a facet. To get it, take a prism over a simplex with two interior points connected by an edge. In the same way, one could  take three interior points in a plane parallel to the top and bottom facets, to get a $5$"=subdivision whose tight span consists of three tetrahedra all sharing an edge. Taking as point configuration the vertices of two simplices, one of them in the interior of the other, one can get a $5$"=subdivision whose tight span consists of four tetrahedra all sharing a vertex. Altogether, we have described all pure three"=dimensional complexes that may occur as the tight span of a $5$"=subdivision.
\end{exmp}

\begin{figure}[htb]\centering
  \includegraphics[width=0.45\textwidth]{5-subdivs-3d.3}
 \includegraphics[width=0.5\textwidth]{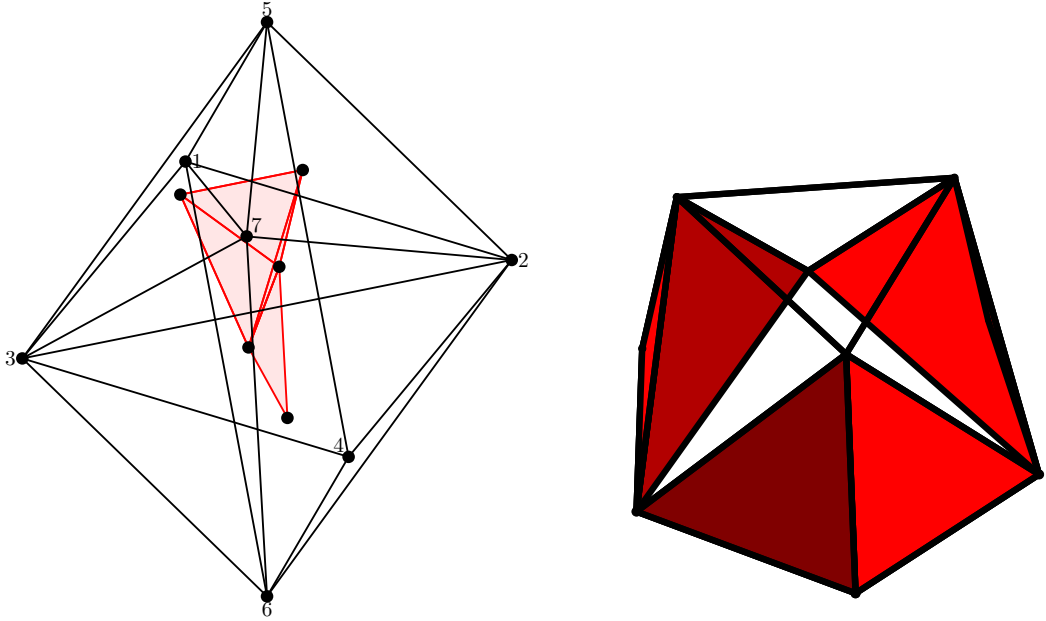}
    \caption{Two $5$"=subdivisions.}
  \label{fig:5-subdiv-comb}
\end{figure}

\begin{exmp}
An example of a subdivision with non"=pure tight span is given in Figure~\ref{fig:5-subdiv-comb} (left). Its tight span is a tetrahedron with a triangle glued at an edge. The point configuration~$\cA$ consists of the six vertices of an octahedron together with an interior point. (Note that the interior point cannot be chosen arbitrarily in this case since one might get a subdivision that is not coarsest.) The subdivision $\subdiv$ of~$\cA$ with maximal faces $\{2,3,4,5,7\}$, $\{1,2,5,7\}$, $\{1,3,5,7\}$, $\{2,3,4,6\}$, and $\{1,2,3,6,7\}$ can be shown to be coarsest and its tight span is as desired, as can be seen from Figure~\ref{fig:5-subdiv-comb}.
Our last example is a $5$"=subdivision with a two"=dimensional tight span that is not planar. In Figure~\ref{fig:5-subdiv-comb} (right), we depicted a polytope subdivided into three simplices and one (rotated) prism over a triangle; this picture was created using \texttt{polymake}~\cite{polymake} and \texttt{JavaView}~\cite{javaview}. Reflecting this complex at the hexagonal facet, one arrives at a polytope with 12 vertices subdivided into six simplices and two triangular prisms. The union of each pair of simplices is convex, hence we can replace them by their union, arriving at a $5$"=subdivision. The tight span of this $5$"=subdivision consists of three triangles that share a common edge.
\end{exmp}

\begin{rem}
\begin{enumerate}
%\item In the proof of Lemma~\ref{lem:5-subdivs}, we have the first case of a simply connected polyhedral complex that cannot occur as a tight span of a $k$"=subdivision and is not excluded by Proposition~\ref{prop:k-no-n-gon} or Proposition~\ref{prop:1-connected}.
\item The examples in Figure~\ref{fig:5-subdiv-2} show that all simply connected  polyhedral complexes with five vertices whose graphs are $2$"=connected and whose maximal faces are all triangles can occur as the tight span of some point configuration. In fact, it can be shown that this true for such complexes with an arbitrary number of vertices.
\item The proof of Lemma~\ref{lem:5-subdivs} can be extended to show that the tight span of any $k$"=subdivision cannot be a $(k-1)$"=gon glued with a triangle.
\end{enumerate}
\end{rem}
As we have seen in Lemma~\ref{lem:4-subdivs} and Example~\ref{ex:4-subdiv-3d}, all three"=dimensional polytopes with up to five vertices can appear as tight spans of $k$"=subdivisions. Since all polytopes can occur as the tight span of some subdivision by Theorem~\ref{thm:ts-all}, it seems natural to ask if all polytopes of dimension three or higher can occur as the tight span of some $k$"=subdivision. The following theorem answers this question negatively.

\begin{thm}
Not all polytopes with dimension three or higher can occur as tight spans of a coarsest subdivision of some point configuration.

Especially, there does not exist a point configuration $\cA$ and a subdivisions $\Sigma$ of~$\cA$ such that the tight span $\tightspan \Sigma \cA$ is a prism over a triangle.
\end{thm}
\begin{proof}
Suppose there exists some point configuration~$\cA$ and a subdivision $\subdiv$ of~$\cA$ such that $\tightspan \subdiv \cA$ is a prism over a triangle. Denote by $F$ the codimension"=three"=cell of~$\subdiv$ corresponding to the prism itself, and by $F_1,F_2,F_3$ the codimension"=one"=cells corresponding to the three parallel edges of $\tightspan \subdiv \cA$. Since $F=F_1\cut F_2 \cut F_3$ is of codimension two in the $F_i$, either $F_1,F_2$, and $F_3$ lie in a common hyperplane $H$, or for each of the hyperplanes  $H_i$ spanned by one of the $F_i$, say $F_1$, the relative interiors of $\conv F_2$ and $\conv F_3$ lie on the same side of~$H_1$. In the first case, the hyperplane $H$ defines a 2"=split of~$\cA$, since the intersection of $H$ with the boundary of $\conv \cA$ equals the intersection of $\conv F_1\union \conv F_2 \union  \conv F_3$ with the boundary and hence cannot produce additional vertices. Obviously, this 2"=split coarsens $\subdiv$.

In the second case, we denote by $H_i^+$ that of the two (closed) halfspaces defined by~$H_i$ that contains the two other faces $F_j$. Obviously, $C:=\conv \cA \cut H_1^+\cut H_2^+\cut H_3^+$ is convex and the union of three maximal cells of $\subdiv$. So we can define a new subdivision $\subdiv'$ of~$\cA$ by replacing these three cells with $C\cut \cA$. Property (SD2) is obviously fulfilled by $\subdiv'$, and, since $F_1,F_2$, and $F_3$ are facets of $C$, (SD3) also holds for $\subdiv'$. Hence $\subdiv'$ is a valid subdivision that coarsens $\subdiv$.

Altogether, $\subdiv$ cannot be a coarsest subdivision of~$\cA$.
\end{proof}

\section{Matroid Subdivisions}

We will now apply our theory of $k$"=splits to a particular class of polytopes, the hypersimplices, more specifically the study of their matroid subdivisions.

We first give the necessary definitions. We abbreviate $[n]:=\{1,2,\dots,n\}$ and
$\binom{[n]}{k}:=\smallSetOf{X\subseteq[n]}{\card X=k}$.  The $k$th \emph{hypersimplex} in
$\RR^n$ is defined as
\[
\Hypersimplex{k}{n} \ := \ \SetOf{x\in[0,1]^n}{\sum_{i=1}^n x_i=k} \ = \ \conv\SetOf{\sum_{i\in A}e_i}{A\in\binom{[n]}{k}} \,;
\]
so it is an $(n-1)$"=dimensional polytope. If~$\cM$ is a matroid on the set $[n]$, then the corresponding \emph{matroid polytope} is the convex hull
of those $0/1$"=vectors in $\RR^n$ which are characteristic functions of the bases of~$\cM$.  For a background on matroids, see the monographs of  White~\cite{White86,White92}. A subdivision $\Sigma$ of $\Hypersimplex kn$ is called a \emph{matroid subdivision} if all $F\in \subdiv$ are matroid polytopes.

A particular example of a matroid is obtained in the following way: Consider a point configuration $\cA\subset\RR^r$. The \emph{matroid of affine dependencies}~$\cM(\cA)$ of $\cA$ is defined by taking as independent sets of~$\cM(\cA)$ the affinely independent subconfigurations of $\cA$. So the bases of~$\cM(\cA)$ are the maximal affinely independent subsets of~$\cA$.

\begin{rem}\label{rem:dressian}
\begin{enumerate}
\item Gel{\cprime}fand, Goresky, MacPherson, and Serganova gave the following characterization of matroid subdivisions~\cite[Theorem~4.1]{GGMS87}: A polytopal subdivision $\subdiv$ of  $\Hypersimplex{k}{n}$ is a matroid subdivision if and only if the $1$"=skeleton of~$\Sigma$ coincides with the $1$"=skeleton of $\Hypersimplex{k}{n}$.
\item The set of all weight functions $w: \cA(\Hypersimplex kn) \to \RR$ that define (regular) matroid subdivisions is the support of a polyhedral fan which is a subfan of the secondary fan of $\Hypersimplex kn$. Speyer~\cite{MR2448909} showed that the set of all those weight vectors is equal to the space of all tropical Plücker vectors, which form the \emph{Dressian} $\Dr(k,n)$. This space includes as a subspace the \emph{tropical Grassmannian} $\Gr(k,n)$ of Speyer and Sturmfels~\cite{SpeyerSturmfels04}, the space of all tropicalized Plücker vectors, or, equivalently, the tropicalization of the usual Grassmannian of all $k$"=dimensional subspaces of an $n$"=dimensional vector space.
\end{enumerate}
\end{rem}

We now recall the description of the 2"=splits of $\Hypersimplex{k}{n}$ given in \cite[Section~5]{MR2502496}. For a triplet $(A,B;\mu)$ with $\emptyset\ne A,B\subsetneq[n]$, $A\cup B=[n]$, $A\cap B=\emptyset$ and $\mu\in\NN$ the hyperplane defined by
\begin{equation}\label{eq:hypersimplex:split}
  \mu \sum_{i\in A} x_i \ = \ (k-\mu) \sum_{i\in B} x_i \,
\end{equation}
is called the \emph{$(A,B;\mu)$"=hyperplane}. Since $\sum_{i=1}^nx_i=k$ for all $x\in \Hypersimplex kn$ this hyperplane can equivalently be described as $\sum_{i\in B}x_i=\mu$.
The 2"=splits of $\Hypersimplex kn$ are now given by all $(A,B,\mu)$"=hyperplanes with $k-\mu+1\le\card A\le n-\mu-1$ and $1\le\mu\le k-1$; see \cite[Lemma~5.1, Proposition~5.2]{MR2502496}. 

\begin{rem}
In \cite[Lemma~7.4]{MR2502496}, it was shown that all 2"=splits of $\Hypersimplex kn$  are matroid subdivisions. So the weight vectors of 2"=splits of $\Hypersimplex kn$ correspond to rays of the Dressian $\Dr(k,n)$. Even more is true: All weight functions in the 2"=split complex of $\Hypersimplex kn$ define matroid subdivisions \cite[Theorem~7.8]{MR2502496}. This gives us the description of a subcomplex of the Dressian. This was used by Jensen, Joswig, Sturmfels, and the author to give a bound on the dimension of the space of all tropical Plücker vectors $\Dr(3,n)$ \cite[Theorem~3.6]{HJJS}.
\end{rem}

We will now construct a class of $3$"=splits of hypersimplices:

\begin{prop}\label{prop:k-splits-from-splits-hs}
Let $A_1\dot\union A_2\dot\union A_3=[n]$ be a partition of $[n]$ into three parts and $\mu_1,\mu_2,\mu_3\in\NN$ such that $\mu_1+\mu_2+\mu_3=k$ and $1\leq\mu_j\leq |A_j|-1$, $j\in\{1,2,3\}$. Then the $([n]\setminus A_j,A_j,\mu_j)$"=hyperplanes define two $3$"=splits of $\Hypersimplex kn$.
\end{prop}
\begin{proof}
We define the polytopes
\begin{align*}
P_1&:=\SetOf{x\in\Hypersimplex kn}{\sum_{i\in A_3}x_i\leq \mu_3\text{ and }\sum_{i\in A_2}x_i\geq \mu_2}\,,\\
P_2&:=\SetOf{x\in\Hypersimplex kn}{\sum_{i\in A_1}x_i\leq \mu_1\text{ and }\sum_{i\in A_3}x_i\geq \mu_3}\,,\\
P_3&:=\SetOf{x\in\Hypersimplex kn}{\sum_{i\in A_2}x_i\leq \mu_2\text{ and }\sum_{i\in A_1}x_i\geq \mu_1}\,,
\end{align*}
each bounded by two of the $([n]\setminus A_j,A_j,\mu_j)$"=hyperplanes. We claim that $P_1,P_2$, and $P_3$ form the maximal cells of a subdivision $\subdiv$ of $\Hypersimplex kn$. Consider some point $x\in\Hypersimplex kn$. Since $\sum_{i=1}^nx_i=k=\mu_1+\mu_2+\mu_3$, there has to be at least one $j\in\{1,2,3\}$ such that $\sum_{i\in A_j}x_i\leq \mu_j$ and one $l\in\{1,2,3\}\setminus\{j\}$ such that $\sum_{i\in A_l}x_i\geq \mu_l$, hence $x$ is in one of $P_1,P_2,P_3$, so (SD2) is fulfiled. Furthermore, $P_1,P_2$ and $P_3$ have distinct relative interiors by definition, so (SD3) is also fulfiled. Finally, the intersection $P_1\cap P_2\cap P_3$ is equal to the $(n-3)$"=dimensional polytope $\SetOf{x\in\Hypersimplex kn}{\sum_{i\in A_j}x_i=\mu_j \text{ for }j\in\{1,2,3\}}$. We deduce that $\subdiv$ is a $3$"=split. A second $3$"=split may be obtained by change each ``$\leq$'' to a ''$\geq$`` and vice versa in the definition of the $P_j$.
\end{proof}

\begin{cor}
The number of $3$"=splits of the hypersimplex $\Hypersimplex kn$ is at least
\begin{equation}\label{eq:n-3-splits}
\frac 13 \sum_{\alpha=2}^{n-4}\sum_{\beta=2}^{n-\alpha-2}\mu^n_k(\alpha,\beta)\binom n{\alpha}\binom{n-\alpha}{\beta}\,,
\end{equation}
where $\mu^n_k(\alpha,\beta)=\sum_{i=1}^{\min(\alpha-1,k-2)}\big(\min(\beta-1,k-i-1)-\max(0,k-i-(n-\alpha-\beta))\big)$.
\end{cor}
\begin{proof}
The number of partitions of $[n]$ into three parts $A_1,A_2,A_3$ where one part has $\alpha$ elements, one has $\beta$ elements, and the last has $n-\alpha-\beta$ elements is $\frac16 \binom n{\alpha}\binom{n-\alpha}{\beta}$. The value $\mu^n_k(\alpha,\beta)$ counts the number of possible choices for $\mu_1,\mu_2,\mu_3$ with $\mu_1+\mu_2+\mu_3=k$ and $1\leq\mu_1\leq \alpha-1$, $1\leq\mu_2\leq\beta-1$, and $1\leq\mu_3\leq n-\alpha-\beta-1$. Now \eqref{eq:n-3-splits} follows by summing over all possibilities for $\alpha$ and $\beta$, taking into account that we need $\alpha\geq2$, $\beta\geq2$, and $n-\alpha-\beta\geq2$.

To compute the value $\mu^n_k(\alpha,\beta)$, we sum over all possible choices for $\mu_1$ and count the so-arising possibilities for $\mu_2$. Since $\mu_2,\mu_3\geq1$, we get $\mu_1\leq k-2$ and, similarly, we get $\mu_2\leq k-\mu_1-1$. To ensure that $\mu_3\leq n-\alpha-\beta-1$, we also need $\mu_2\geq k-\mu_1-(n-\alpha-\beta)+1$. This shows the formula.
\end{proof}

For $k=3$ (and $n\geq 6$; there do not exist any $3$"=splits for $\Hypersimplex kn$ with $n\leq 6$), we obviously have $\mu^n_3(\alpha,\beta)=1$ and for $k=4$ we get the simpler formula $\mu^n_4(\alpha,\beta)=3-\card{\smallSetOf{j\in\{\alpha,\beta,n-\alpha-\beta\}}{j=2}}$.

\begin{thm}\label{thm:k-splits-matroid-subdivs}
The $3$"=splits of $\Hypersimplex kn$ constructed in Proposition~\ref{prop:k-splits-from-splits-hs} are matroid subdivisions.
\end{thm}

For the proof we need the following notions from linear algebra. Let $V$ be vector space. A point configuration $\cA\subset V$ is said to be in \emph{general position} if any  $S\subset \cA$ with $\card S\leq \dim V+1$ is affinely independent. A family $\smallSetOf{\cA_i}{i\in I}$ of point configurations in $V$ is said to be in \emph{relative general position} if for each affinely dependent set $S\subseteq \bigcup_{i\in I} \cA_i$ with $\card S\leq \dim V+1$ there exists some $i\in I$ such that $S\cap \cA_i$ is affinely dependent.
We furthermore need the following result \cite[Lemma~7.3]{MR2502496}: Let $\cA\subset \RR^{k-1}$ be a point configuration such that there exists a family $\smallSetOf{\cA_i}{i\in I}$ of point configurations in relative general position such that each $\cA_i$ is in general position as a subset of
  $\aff \cA_i$ and such that  $\cA=\bigcup_{i\in I} \cA_i$. Then the set of bases of $\cM(\cA)$ is given by
  \begin{align}\label{eq:matroid-general-position}
    \SetOf{B\subset \cA}{\card B=k\text{ and } |(B\cap \cA_i)|\leq \dim (\aff \cA_i)+1\text{ for all }i\in I} \, .
  \end{align}

\begin{proof}[Proof of Theorem~\ref{thm:k-splits-matroid-subdivs}]
%Consider the partition $\dot\bigcup_{j=1}^l A_j=[n]$ of $[n]$ and $\mu_j$ ($1\leq j\leq n$) satisfying the conditions of Proposition~\ref{prop:k-splits-from-splits-hs}. Define a point configuration~$\cA$ in $\RR^{l-1}$ by choosing $l$ affinely independent points $a_1,\dots a_l$ and taking $a_i$ $\mu_i$ times. The bases of the matroid of affine dependences of~$\cA$ are in bijection with the vertices of
Each full"=dimensional face $F$ of a subdivision obtained by the construction in Proposition~\ref{prop:k-splits-from-splits-hs} is the intersection of~$\cA(\Hypersimplex kn)$ with some $H_j^+$ where $H_i$ is the $(A_j,[n]\setminus A_j,\mu_j)$"=hyperplane. So, without loss of generality, let $F=\cA(\Hypersimplex kn) \cut \Cut_{j\in J} H_j^+$ for some $J\subset \{1,2,3\}$. The elements of $F$ are all $0/1$"=vectors $x$ of length $n$ with $k$ ones that fulfill $\sum_{i\in A_j} x_i\leq \mu_i$ for all $j\in J$. We will construct a point configuration~$\cA\subset\RR^{k-1}$ with $n$ points such that $\conv F$ is the matroid polytope~$\cM(\cA)$.

For each $j\in J$ we choose a $(\mu_j-1)$"=dimensional affine subspace $U_j$ of $\RR^{k-1}$ such that $U_i\cut U_j=\emptyset$ for all $i,j\in J$. This is possible since $\sum_{j\in J}\mu_j\leq \sum_{j=1}^l \mu_j=k$. Now we choose for each $j\in J$ a point configuration~$\cA_j\subset U_j$ with $\card{A_j}$ points such that~$\cA_j$ is in general position in $U_i$ and such that the family $\smallSetOf{\cA_j}{j\in J}$ is in relative general position. The final $n-\card{{\dot\union_{j\in J} A_j}}$ points of~$\cA$ are chosen in general position in $\RR^{k-1}$. By the discussion above, the bases of~$\cM(\cA)$ are those $k$"=element subsets of~$\cA$ whose intersection with~$\cA_j$ has cardinality smaller or equal to $\mu_j$ for all $j\in J$. This shows the claim.
\end{proof}

\begin{rem}
Together with the construction in Proposition~\ref{prop:k-splits-from-splits-hs}, Theorem~\ref{thm:k-splits-matroid-subdivs} gives us a lot of new rays for the Dressian $\Dr(k,n)$ (whose weight vectors can be constructed as in the proof of Theorem~\ref{thm:k-regular}). This is a further step in the understanding of this space of tropical Plücker vectors. Via the complete computation of $\Dr(3,6)$~\cite{SpeyerSturmfels04} and $\Dr(3,7)$~\cite{HJJS}, we see that these are not all rays, even if $k=3$; but this gives us at least some more information about the Dressian in the general case.
\end{rem}

\section{Open Questions}

We have discussed some conditions on when a polyhedral complex can be the tight span of some $k$"=subdivision. However, we also gave examples that these conditions are not sufficient. For complexes with a sole maximal cell, we showed that the only possibility in dimension two is a triangle, and that in dimension three not all polytopes may occur. This naturally leads to the following question.

\begin{qst}
  Which polyhedral complexes, especially, which polytopes occur as tight spans of $k$"=subdivisions?
\end{qst}

Especially, it might be interesting to define and analyze special classes of $k$"=subdivisions other than $k$"=splits. 

\begin{qst}
  Which polytopes are totally $k$"=splittable?% Which polytopes are totally semi"=splittable?
\end{qst}

The answer for this question might lead to interesting new classes of polytopes, the class of all totally $3$"=splittable polytopes, all totally $4$"=splittable polytopes, and so on. This would help to get new insights into the structure of secondary polytopes. Especially, since for the class of totally $2$"=splittable polytopes all secondary polytopes are known, a classification of totally $k$"=splittable polytopes for small $k\geq 3$ could lead to explicit computations of some secondary polytopes.

In \cite{MR2502496} it was shown that the 2"=split complex of $\Hypersimplex kn$ is a subcomplex of the complex of all matroid subdivisions of $\Hypersimplex kn$. As $3$"=splits are also matroid subdivisions, the following seams natural to ask:

\begin{qst}
Are refinements of $3$"=splits (or $l$"=splits) of $\Hypersimplex kn$ again matroid subdivisions? 
\end{qst}

\bibliographystyle{amsplain}
\bibliography{diss,main}

\end{document}